\documentclass{elsart}
\usepackage{amsfonts}

\usepackage{graphicx}
\usepackage{pstricks}
\usepackage{amsmath}
\usepackage{amsmath}
\usepackage{amsxtra}
\usepackage{graphicx}
\usepackage{pstricks}
\usepackage{amsmath}
\usepackage{amsmath}
\usepackage{amsxtra}

\newtheorem{theorem}{Theorem}[section]
\newtheorem{lemma}[theorem]{Lemma}
\newtheorem{proposition}[theorem]{Proposition}
\newtheorem{corollary}[theorem]{Corollary}

\newtheorem{conjecture}[theorem]{Conjecture}
\newtheorem{definition}[theorem]{Definition}

\newtheorem{example}[theorem]{Example}
\newtheorem{remark}[theorem]{Remark}

\numberwithin{equation}{section}
\begin{document}

\begin{frontmatter}



\title{Hyponormality and Subnormality for Powers
of Commuting Pairs of Subnormal Operators}


\author{Ra\'{u}l E. Curto}
\address{Department of Mathematics, The University of Iowa, Iowa City, Iowa
52242}
\ead{rcurto@math.uiowa.edu}
\ead[url]{http://www.math.uiowa.edu/\symbol{126}rcurto/}
\author{Sang Hoon Lee}
\address{Department of Mathematics, The University of Iowa, Iowa City, Iowa
52242}
\ead{shlee@math.skku.ac.kr}
%\urladdr{}
\author{Jasang Yoon}
\address{Department of Mathematics, Iowa State University, Ames, Iowa 50011}
\ead{jyoon@iastate.edu}
\ead[url]{http://www.public.iastate.edu/\symbol{126}jyoon/}
\thanks{Research partially supported by NSF Grants DMS-0099357 and DMS-0422952.}
%\MSC{Primary 47B20, 47B37, 47A13, 28A50; Secondary 44A60, 47-04, 47A20}
%\keywords

\begin{abstract}
Let $\mathfrak{H}_{0}$ (resp. $\mathfrak{H}_{\infty }$) denote the class of
commuting pairs of subnormal operators on Hilbert space (resp. subnormal
pairs), and for an integer $k\geq 1$ let $\mathfrak{H}_{k}$ denote the class
of $k$-hyponormal pairs in $\mathfrak{H}_{0}$. \ We study the hyponormality
and subnormality of powers of pairs in $\mathfrak{H}_{k}$. \ We first show
that if $(T_{1},T_{2})\in \mathfrak{H}_{1}$, the pair $(T_{1}^{2},T_{2})$
may fail to be in $\mathfrak{H}_{1}$. \ Conversely, we find a pair 
$(T_{1},T_{2})\in \mathfrak{H}_{0}$ such that 
$(T_{1}^{2},T_{2})\in \mathfrak{H}_{1}$ but $(T_{1},T_{2})\notin \mathfrak{H}_{1}$.
 \ Next, we show that there exists a
pair $(T_{1},T_{2})\in \mathfrak{H}_{1}$ such that $T_{1}^{m}T_{2}^{n}$ is
subnormal (all $m,n\geq 1$), but $(T_{1},T_{2})$ is not in $\mathfrak{H}%
_{\infty }$; this further stretches the gap between the classes $\mathfrak{H}_{1}$ and 
$\mathfrak{H}_{\infty}$.
\ Finally, we prove that there exists a large class of $2$-variable weighted shifts 
$(T_{1},T_{2})$ (namely those pairs in $\mathfrak{H}_{0}$ whose cores are of tensor 
form (cf. Definition \ref{tc})), for which the subnormality of $(T_{1}^{2},T_{2})$
and $(T_{1},T_{2}^{2})$ does imply the subnormality of $(T_{1},T_{2})$.
\end{abstract}

\begin{keyword}
jointly hyponormal pairs \sep subnormal pairs \sep $2$-variable weighted
shifts \sep powers of commuting pairs of subnormal operators

\textit{2000 Mathematics Subject Classification.}  Primary 47B20, 47B37, 47A13, 28A50; 
Secondary 44A60, 47-04, 47A20 
%\PACS 
\end{keyword}
%\maketitle
\end{frontmatter}% main text

\section{Introduction}

\label{Int}

The Lifting Problem for Commuting Subnormals (LPCS) asks for necessary and
sufficient conditions for a pair of subnormal operators on Hilbert space to
admit commuting normal extensions. \ It is well known that the commutativity
of the pair is necessary but not sufficient (\cite{Abr}, \cite{Lu1}, \cite%
{Lu2}, \cite{Lu3}), and it has recently been shown that the joint
hyponormality of the pair is necessary but not sufficient \cite{CuYo1}, thus
disproving the conjecture in \cite{CMX}. \ An abstract answer to the Lifting
Problem was obtained in \cite{CLY}, by stating and proving a multivariable
analogue of the Bram-Halmos criterion for subnormality, and then showing
concretely that no matter how $k$-hyponormal a pair might be, it may still
fail to be subnormal. \ While this provides new insights into the LPCS, it
stops short of identifying other types of conditions that, together with
joint hyponormality, may imply subnormality.

Our previous work (\cite{CLY}, \cite{CuYo1}, \cite{CuYo2}, \cite{CuYo3}, %
\cite{Yoo}) has revealed that the nontrivial aspects of the LPCS are best
detected within the class $\mathfrak{H}_{1}$ of commuting hyponormal pairs
of subnormal operators; we thus focus our attention on this class. \ More
generally, we will denote the class of commuting pairs of subnormal
operators on Hilbert space by $\mathfrak{H}_{0}$, the class of subnormal
pairs by $\mathfrak{H}_{\infty }$, and for an integer $k\geq 1$ the class of 
$k$-hyponormal pairs in $\mathfrak{H}_{0}$ by $\mathfrak{H}_{k}$. \ Clearly, 
$\mathfrak{H}_{\infty }\subseteq ...\subseteq \mathfrak{H}_{k}\subseteq
...\subseteq \mathfrak{H}_{2}\subseteq \mathfrak{H}_{1}\subseteq \mathfrak{H}%
_{0}$; the main results in \cite{CuYo1} and \cite{CLY} show that these
inclusions are all proper. \ (The LPCS thus asks for necessary and
sufficient conditions for a pair $\mathbf{T\in }\mathfrak{H}_{0}$ to be in $%
\mathfrak{H}_{\infty }$.) \ 

In \cite{Fra}, E. Franks proved that if $\mathbf{T}\equiv (T_{1},T_{2})\in 
\mathfrak{H}_{0}$ and $p(\mathbf{T})$ is subnormal for all polynomials $p\in 
\mathbb{C}[z]$ with $\deg p\leq 5$, then $\mathbf{T}$ is necessarily
subnormal. \ Motivated in part by this result, and in part by J. Stampfli's
work in \cite{Sta} and \cite{Sta2}, in this article we consider the role of
the powers of a pair in ascertaining its subnormality. \ Clearly, if $%
\mathbf{T}\equiv (T_{1},T_{2})\in \mathfrak{H}_{\infty }$, and if $m,n\geq 1$%
, then $\mathbf{T}^{(m,n)}:=(T_{1}^{m},T_{2}^{n})\in \mathfrak{H}_{\infty }$%
, and therefore $T_{1}^{m}T_{2}^{n}$ is a subnormal operator. \ It is thus
natural to ask whether the subnormality of both $\mathbf{T}^{(2,1)}$ and $%
\mathbf{T}^{(1,2)}$ can force the subnormality of $\mathbf{T}$.

Our first main result shows that the class $\mathfrak{H}_{1}$ is not
invariant under squares, as follows: we construct a pair $\mathbf{T}\equiv
(T_{1},T_{2})\in \mathfrak{H}_{1}$ such that $\mathbf{T}%
^{(2,1)}=(T_{1}^{2},T_{2})\notin \mathfrak{H}_{1}$ (Theorem \ref{firstmain}%
). \ Conversely, we find a pair $\mathbf{T\in }\mathfrak{H}_{0}$ such that $%
\mathbf{T}^{(2,1)}=(T_{1}^{2},T_{2})\in \mathfrak{H}_{1}$ but $\mathbf{T}%
\notin \mathfrak{H}_{1}$. \ We then show that for a large class of commuting
pairs of subnormal operators, the subnormality of both $\mathbf{T}^{(2,1)}$
and $\mathbf{T}^{(1,2)}$ does force the subnormality of $\mathbf{T}$. \
Concretely, if $\mathbf{T}\in \mathcal{TC}$, the class of all $2$-variable
weighted shifts $\mathbf{T}\in \mathfrak{H}_{0}$ whose cores are of \textit{%
tensor form }(see Definition \ref{tc} below), then $\mathbf{T}^{(1,2)}\in 
\mathfrak{H}_{\infty }\Leftrightarrow \mathbf{T}^{(2,1)}\in \mathfrak{H}%
_{\infty }\Leftrightarrow \mathbf{T}\in \mathfrak{H}_{\infty }$ (Theorem \ref%
{Thm1}). \ Our results thus seem to indicate that the subnormality of $%
\mathbf{T}^{(2,1)},\mathbf{T}^{(1,2)}$ may very well be essential in
determining the subnormality of $\mathbf{T}$ within the class $\mathfrak{H}%
_{0}$ (Conjecture \ref{conjecture}). \ Next, we prove that it is possible
for a pair $\mathbf{T}\in \mathfrak{H}_{1}$ to have all powers $%
T_{1}^{m}T_{2}^{n}\;(m,n\geq 1)$ subnormal, without being subnormal (Example %
\ref{four}). \ This provides further evidence that the gap between the
classes $\mathfrak{H}_{\infty }$ and $\mathfrak{H}_{1}$ is fairly large.

To prove our results, we resort to tools introduced in previous work (e.g.,
the Six-point Test to check hyponormality (Lemma \ref{joint hypo}) and the
Backward Extension Theorem for $2$-variable weighted shifts (Lemma \ref%
{backext})), together with a new direct sum decomposition for powers of $2$%
-variable weighted shifts which parallels the decomposition used in \cite%
{CuP} to analyze $k$-hyponormality for powers of (one-variable) weighted
shifts. \ Specifically, we split the ambient space $\ell ^{2}(\mathbb{Z}%
_{+}^{2})$ as an orthogonal direct sum $\mathcal{H}^{0}\oplus \mathcal{H}%
^{1} $, where $\mathcal{H}^{m}:=\bigvee_{k=0}^{\infty
}\{e_{(j,2k+m)}:j=0,1,2,\cdots \}\;(m=0,1)$. \ Each of the subspaces $%
\mathcal{H}^{0}$ and $\mathcal{H}^{1}$ reduces $T_{1}$ and $T_{2}$, and $%
\mathbf{T}^{(1,2)}$ is subnormal if and only if each of $\mathbf{T}%
^{(1,2)}|_{\mathcal{H}^{0}}$ and $\mathbf{T}^{(1,2)}|_{\mathcal{H}^{1}}$ is
subnormal (cf. Figure \ref{double}).

We devote the rest of this section to establishing our basic terminology and
notation. \ Let $\mathcal{H}$ be a complex Hilbert space and let $\mathcal{B}%
(\mathcal{H})$ denote the algebra of bounded linear operators on $\mathcal{H}
$. We say that $T\in \mathcal{B}(\mathcal{H})$ is \textit{normal} if $%
T^{\ast }T=TT^{\ast },$ \textit{subnormal} if $T=N|_{\mathcal{H}}$, where $N$
is normal and $N(\mathcal{H})\mathcal{\subseteq H}$, and \textit{hyponormal}
if $T^{\ast }T\geq TT^{\ast }$. \ For $S,T\in \mathcal{B}(\mathcal{H})$ let $%
[S,T]:=ST-TS$.\ \ We say that an $n$-tuple $\mathbf{T}\equiv (T_{1},\cdots
,T_{n})$ of operators on $\mathcal{H}$ is (jointly) \textit{hyponormal} if
the operator matrix 
\begin{equation*}
\lbrack \mathbf{T}^{\ast },\mathbf{T]:=}\left( 
\begin{array}{llll}
\lbrack T_{1}^{\ast },T_{1}] & [T_{2}^{\ast },T_{1}] & \cdots & [T_{n}^{\ast
},T_{1}] \\ 
\lbrack T_{1}^{\ast },T_{2}] & [T_{2}^{\ast },T_{2}] & \cdots & [T_{n}^{\ast
},T_{2}] \\ 
\text{ \thinspace \thinspace \quad }\vdots & \text{ \thinspace \thinspace
\quad }\vdots & \ddots & \text{ \thinspace \thinspace \quad }\vdots \\ 
\lbrack T_{1}^{\ast },T_{n}] & [T_{2}^{\ast },T_{n}] & \cdots & [T_{n}^{\ast
},T_{n}]%
\end{array}%
\right)
\end{equation*}%
is positive on the direct sum of $n$ copies of $\mathcal{H}$ (cf. \cite{Ath}%
, \cite{CuLe1}, \cite{CMX}). \ The $n$-tuple $\mathbf{T}$ is said to be 
\textit{normal} if $\mathbf{T}$ is commuting and each $T_{i}$ is normal, and 
$\mathbf{T}$ is \textit{subnormal }if $\mathbf{T}$ is the restriction of a
normal $n$-tuple to a common invariant subspace. \ Finally, we say that a
pair $\mathbf{T}\equiv (T_{1},T_{2})$ is $2$-hyponormal if $\mathbf{T}$ is
commuting and $(T_{1},T_{2},T_{1}^{2},T_{1}T_{2},T_{2}^{2})$ is hyponormal.
\ Clearly, normal $\Rightarrow $ subnormal $\Rightarrow $ $2$-hyponormal $%
\Rightarrow $ hyponormal.

The Bram-Halmos criterion for subnormality states that an operator $T\in 
\mathcal{B}(\mathcal{H})$ is subnormal if and only if 
\begin{equation*}
\sum_{i,j}(T^{i}x_{j},T^{j}x_{i})\geq 0
\end{equation*}%
for all finite collections $x_{0},x_{1},\cdots ,x_{k}\in \mathcal{H}$ (\cite%
{Bra}, \cite{Con}). \ Using Choleski's algorithm for operator matrices, it
is easy to verify that this condition is equivalent to the assertion that
the $k$-tuple $(T,T^{2},\cdots ,T^{k})$ is hyponormal for all $k\geq 1$.

For $\alpha \equiv \{\alpha _{n}\}_{n=0}^{\infty }$ a bounded sequence of
positive real numbers (called \textit{weights}) let $W_{\alpha }:\ell ^{2}(%
\mathbb{Z}_{+})\rightarrow \ell ^{2}(\mathbb{Z}_{+})$ be the associated
unilateral weighted shift, defined by $W_{\alpha }e_{n}:=\alpha
_{n}e_{n+1}\;($all $n\geq 0)$, where $\{e_{n}\}_{n=0}^{\infty }$ is the
canonical orthonormal basis in $\ell ^{2}(\mathbb{Z}_{+}).$ \ For notational
convenience, we will often write $shift(\alpha _{0},\alpha _{1},\cdots )$ to
denote $W_{\alpha }$. \ In particular, we shall let $U_{+}:=shift(1,1,\cdots
)$ ($U_{+}$ is the (unweighted) unilateral shift) and $S_{a}:=shift(a,1,1,%
\cdots ).$ \ For a weighted shift $W_{\alpha }$, \textit{the moments of }$%
\alpha $ are given by 
\begin{equation*}
\gamma _{k}\equiv \gamma _{k}(\alpha ):=%
\begin{cases}
1 & \text{if }k=0 \\ 
\alpha _{0}^{2}\cdots \alpha _{k-1}^{2} & \text{if }k>0.%
\end{cases}%
\end{equation*}%
It is easy to see that $W_{\alpha }$ is never normal, and that it is
hyponormal if and only if $\alpha _{0}\leq \alpha _{1}\leq \cdots $. \
Similarly, consider double-indexed positive bounded sequences $\alpha _{%
\mathbf{k}},\beta _{\mathbf{k}}\in \ell ^{\infty }(\mathbb{Z}_{+}^{2})$, $%
\mathbf{k}\equiv (k_{1},k_{2})\in \mathbb{Z}_{+}^{2}:=\mathbb{Z}_{+}\times 
\mathbb{Z}_{+}$, and let $\ell ^{2}(\mathbb{Z}_{+}^{2})$ be the Hilbert
space of square-summable complex sequences indexed by $\mathbb{Z}_{+}^{2}$.
\ (Recall that $\ell ^{2}(\mathbb{Z}_{+}^{2})$ is canonically isometrically
isomorphic to $\ell ^{2}(\mathbb{Z}_{+})\bigotimes \ell ^{2}(\mathbb{Z}_{+})$%
.) \ We define the $2$-variable weighted shift $\mathbf{T}\equiv
(T_{1},T_{2})$ by 
\begin{equation*}
\begin{cases}
T_{1}e_{\mathbf{k}}:=\alpha _{\mathbf{k}}e_{\mathbf{k+}\varepsilon _{1}} \\ 
T_{2}e_{\mathbf{k}}:=\beta _{\mathbf{k}}e_{\mathbf{k+}\varepsilon _{2}},%
\end{cases}%
\end{equation*}%
where $\mathbf{\varepsilon }_{1}:=(1,0)$ and $\mathbf{\varepsilon }%
_{2}:=(0,1)$. \ Clearly, 
\begin{equation}
T_{1}T_{2}=T_{2}T_{1}\Longleftrightarrow \beta _{\mathbf{k+}\varepsilon
_{1}}\alpha _{\mathbf{k}}=\alpha _{\mathbf{k+}\varepsilon _{2}}\beta _{%
\mathbf{k}}\;\;(\text{all }\mathbf{k}\in \mathbb{Z}_{+}^{2}).
\label{commuting}
\end{equation}%
In an entirely similar way one can define multivariable weighted shifts. \ 

A $2$-variable weighted shift $\mathbf{T}\equiv (T_{1},T_{2})$ is called 
\textit{horizontally flat} if $\alpha _{(k_{1},k_{2})}=\alpha _{(1,1)}$ for
all $k_{1},k_{2}\geq 1$; $\mathbf{T}$ \ is called \textit{vertically flat}
if $\beta _{(k_{1},k_{2})}=\beta _{(1,1)}$ for all $k_{1},k_{2}\geq 1$. \ If 
$\mathbf{T}$ is horizontally and vertically flat, then $\mathbf{T}$ is
simply called \textit{flat}.

For an arbitrary $2$-variable weighted shift $\mathbf{T}$, we shall let $%
\mathcal{R}_{ij}(\mathbf{T})$ denote the restriction of $\mathbf{T}$ to $%
\mathcal{M}_{i}\cap \mathcal{N}_{j}$, where $\mathcal{M}_{i}$ (resp. $%
\mathcal{N}_{j})$ is the subspace of $\ell ^{2}(\mathbb{Z}_{+}^{2})$ spanned
by the canonical orthonormal basis vectors associated to indices $\mathbf{k}%
=(k_{1},k_{2})$ with $k_{1}\geq 0$ and $k_{2}\geq i$ (resp. $k_{1}\geq j$
and $k_{2}\geq 0$).

Trivially, a pair of unilateral weighted shifts $W_{\alpha }$ and $W_{\beta
} $ gives rise to a $2$-variable weighted shift $\mathbf{T}\equiv
(T_{1},T_{2}) $, if we let $\alpha _{(k_{1},k_{2})}:=\alpha _{k_{1}}$ and $%
\beta _{(k_{1},k_{2})}:=\beta _{k_{2}}\;$(all $k_{1},k_{2}\in \mathbb{Z}%
_{+}^{2}$ ). \ In this case, $\mathbf{T}$ is subnormal (resp. hyponormal) if
and only if so are $T_{1}$ and $T_{2}$; in fact, under the canonical
identification of $\ell ^{2}(\mathbb{Z}_{+}^{2})$ and $\ell ^{2}(\mathbb{Z}%
_{+})\bigotimes \ell ^{2}(\mathbb{Z}_{+})$, $T_{1}\cong I\bigotimes
W_{\alpha }$ and $T_{2}\cong W_{\beta }\bigotimes I$, and $\mathbf{T}$ is
also doubly commuting. \ For this reason, we do not focus attention on
shifts of this type, and use them only when the above mentioned triviality
is desirable or needed. \ Given $\mathbf{k}\in \mathbb{Z}_{+}^{2}$, the
moment of $(\alpha ,\beta )$ of order $\mathbf{k}$ is 
\begin{equation*}
\gamma _{\mathbf{k}}\equiv \gamma _{\mathbf{k}}(\alpha ,\beta ):=%
\begin{cases}
1 & \text{if }\mathbf{k}=0 \\ 
\alpha _{(0,0)}^{2}\cdots \alpha _{(k_{1}-1,0)}^{2} & \text{if }k_{1}\geq 1%
\text{ and }k_{2}=0 \\ 
\beta _{(0,0)}^{2}\cdots \beta _{(0,k_{2}-1)}^{2} & \text{if }k_{1}=0\text{
and }k_{2}\geq 1 \\ 
\alpha _{(0,0)}^{2}\cdots \alpha _{(k_{1}-1,0)}^{2}\beta
_{(k_{1},0)}^{2}\cdots \beta _{(k_{1},k_{2}-1)}^{2} & \text{if }k_{1}\geq 1%
\text{ and }k_{2}\geq 1.%
\end{cases}%
\end{equation*}%
(We remark that, due to the commutativity condition (\ref{commuting}), $%
\gamma _{\mathbf{k}}$ can be computed using any nondecreasing path from $%
(0,0)$ to $(k_{1},k_{2})$.) \ We now recall a well known characterization of
subnormality for multivariable weighted shifts \cite{JeLu}, due to C. Berger
(cf. \cite[III.8.16]{Con}) and independently established by Gellar and
Wallen \cite{GeWa}) in the single variable case: $\ \mathbf{T\equiv (}%
T_{1},T_{2})$ admits a commuting normal extension if and only if there is a
probability measure $\mu $ (which we call the Berger measure of $\mathbf{T}$%
) defined on the $2$-dimensional rectangle $R=[0,a_{1}]\times \lbrack
0,a_{2}]$ (where $a_{i}:=\left\| T_{i}\right\| ^{2}$) such that $\gamma _{%
\mathbf{k}}=\int_{R}s^{k_{1}}t^{k_{2}}d\mu (s,t),$ for all $\mathbf{k}\in 
\mathbb{Z}_{+}^{2}$. \ In the single variable case, if $W_{\alpha }$ is
subnormal with Berger measure $\xi _{\alpha }$ and $h\geq 1$, and if we let $%
\mathcal{L}_{h}:=\bigvee \{e_{n}:n\geq h\}$ denote the invariant subspace
obtained by removing the first $h$ vectors in the canonical orthonormal
basis of $\ell ^{2}(\mathbb{Z}_{+})$, then the Berger measure of $W_{\alpha
}|_{\mathcal{L}_{h}}$ is $\frac{s^{h}}{\gamma _{h}}d\xi (s)$; alternatively,
if $S:\ell ^{\infty }(\mathbb{Z}_{+})\rightarrow \ell ^{\infty }(\mathbb{Z}%
_{+})$ is defined by 
\begin{equation}
S(\alpha )(n):=\alpha (n+1)\;(\alpha \in \ell ^{\infty }(\mathbb{Z}%
_{+}),n\geq 0),  \label{sa}
\end{equation}%
then 
\begin{equation}
d\xi _{S(\alpha )}(s)=\frac{s}{\alpha _{0}^{2}}d\xi (s).  \label{sam}
\end{equation}

\textit{Acknowledgments}. \ The authors are very grateful to the referee for
several suggestions which helped improved the presentation. \ Most of the examples, 
and some of the proofs in this paper were obtained using calculations with the 
software tool \textit{Mathematica \cite{Wol}.}

\section{The Class $\mathfrak{H}_{1}$ Is Not Invariant Under Squares}

\label{powers}

For a general operator $T$ on Hilbert space, it is well known that the
hyponormality of $T$ does not imply the hyponormality of $T^{2}$ (\cite{Hal}%
). \ However, for a unilateral weighted shift $W_{\alpha }$, the
hyponormality of $W_{\alpha }$ (detected by the condition $\alpha _{k}\leq
\alpha _{k+1}$ for all $k\geq 0$) clearly implies the hyponormality of every
power $W_{\alpha }^{m}\;(m\geq 1)$. $\ $For $2$-variable weighted shifts,
one is thus tempted to expect that a similar result would hold, especially
if we restrict attention to the class $\mathfrak{H}_{1}$ of commuting
hyponormal pairs of subnormal operators. \ Somewhat surprisingly, it is
actually possible to build a $2$-variable weighted shift $\mathbf{T}\in 
\mathbf{\mathfrak{H}_{1}}$ such that $\mathbf{T}^{(2,1)}\notin \mathfrak{H}%
_{1}$, and we do this in this section. \ 

We begin with some basic results. \ First, we recall a hyponormality
criterion for $2$-variable weighted shifts.

\begin{lemma}
$($\cite{bridge}$)$\label{joint hypo} (Six-point Test) \ Let $\mathbf{%
T\equiv (}T_{1},T_{2})$ be a $2$-variable weighted shift, with weight
sequences $\alpha $ and $\beta $. \ Then $\mathbf{T}$ is hyponormal if and
only if 
\begin{equation*}
H_{\mathbf{T}}(\mathbf{k}):=\left( 
\begin{array}{cc}
\alpha _{\mathbf{k}+\mathbf{\varepsilon }_{1}}^{2}-\alpha _{\mathbf{k}}^{2}
& \alpha _{\mathbf{k}+\mathbf{\varepsilon }_{2}}\beta _{\mathbf{k}+\mathbf{%
\varepsilon }_{1}}-\alpha _{\mathbf{k}}\beta _{\mathbf{k}} \\ 
\alpha _{\mathbf{k}+\mathbf{\varepsilon }_{2}}\beta _{\mathbf{k}+\mathbf{%
\varepsilon }_{1}}-\alpha _{\mathbf{k}}\beta _{\mathbf{k}} & \beta _{\mathbf{%
k}+\mathbf{\varepsilon }_{2}}^{2}-\beta _{\mathbf{k}}^{2}%
\end{array}%
\right) \geq 0\quad \text{(all }\mathbf{k}\equiv (k_{1},k_{2})\in \mathbb{Z}%
_{+}^{2}\text{)}.
\end{equation*}
\end{lemma}

\setlength{\unitlength}{1mm} \psset{unit=1mm}

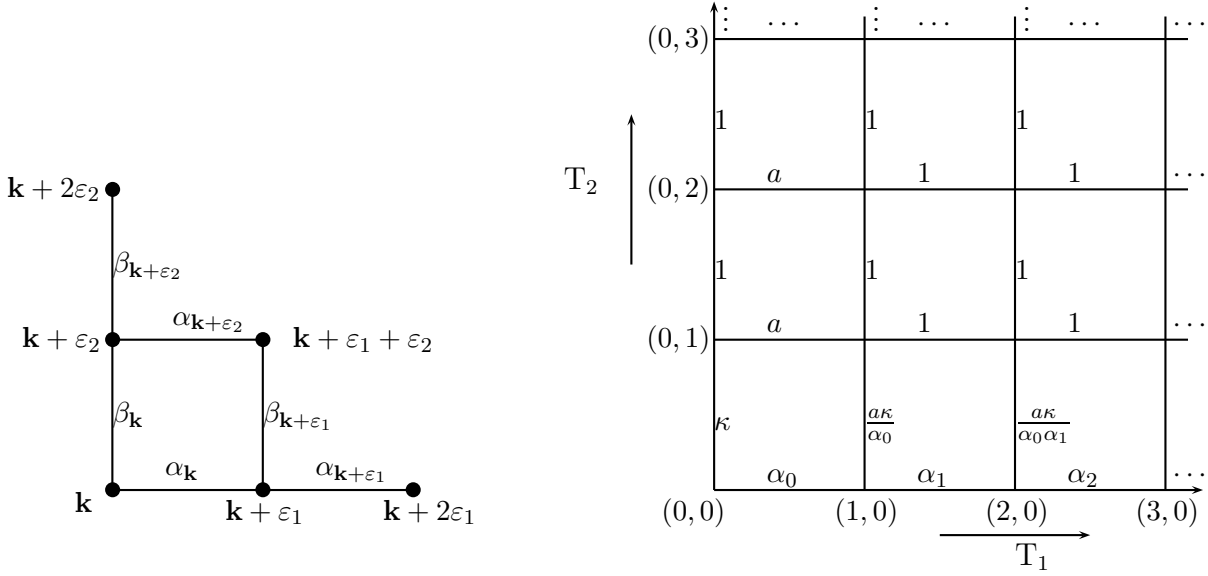
\begin{figure}[th]
\begin{center}
\begin{picture}(160,75)

\psline(10,10)(50,10)
\psline(10,30)(30,30)
\psline(10,10)(10,50)
\psline(30,10)(30,30)

\put(10,10){\pscircle*(0,0){1}} \put(30,10){\pscircle*(0,0){1}}
\put(50,10){\pscircle*(0,0){1}} \put(10,30){\pscircle*(0,0){1}}
\put(30,30){\pscircle*(0,0){1}} \put(10,50){\pscircle*(0,0){1}}

\put(5,7){\footnotesize{$\mathbf{k}$}}
\put(25,6){\footnotesize{$\mathbf{k}+\mathbf{\varepsilon}_{1}$}}
\put(46,6){\footnotesize{$\mathbf{k}+2\mathbf{\varepsilon}_{1}$}}

\put(17,12){\footnotesize{$\alpha_{\mathbf{k}}$}}
\put(37,12){\footnotesize{$\alpha_{\mathbf{k}+\mathbf{\varepsilon}_{1}}$}}

\put(18,32){\footnotesize{$\alpha_{\mathbf{k}+\mathbf{\varepsilon}_{2}}$}}

\put(-2,29){\footnotesize{$\mathbf{k}+\mathbf{\varepsilon}_{2}$}}
\put(-4,49){\footnotesize{$\mathbf{k}+2\mathbf{\varepsilon}_{2}$}}
\put(34,29){\footnotesize{$\mathbf{k}+\mathbf{\varepsilon}_{1}+\mathbf{\varepsilon}_{2}$}}

\put(10,19){\footnotesize{$\beta_{\mathbf{k}}$}}
\put(10,39){\footnotesize{$\beta_{\mathbf{k}+\mathbf{\varepsilon}_{2}}$}}

\put(30,19){\footnotesize{$\beta_{\mathbf{k}+\mathbf{\varepsilon}_{1}}$}}

%\\\\\\

\psline{->}(120,4)(140,4)
\put(130,0){$\rm{T}_1$}
\psline{->}(79,40)(79,60)
\put(70,50){$\rm{T}_2$}

\psline{->}(90,10)(155,10)
\psline(90,30)(153,30)
\psline(90,50)(153,50)
\psline(90,70)(153,70)

\psline{->}(90,10)(90,75)
\psline(110,10)(110,73)
\psline(130,10)(130,73)
\psline(150,10)(150,73)

\put(83,6){\footnotesize{$(0,0)$}}
\put(106,6){\footnotesize{$(1,0)$}}
\put(126,6){\footnotesize{$(2,0)$}}
\put(146,6){\footnotesize{$(3,0)$}}
\put(81,29){\footnotesize{$(0,1)$}}
\put(81,49){\footnotesize{$(0,2)$}}
\put(81,69){\footnotesize{$(0,3)$}}

\put(97,11){\footnotesize{$\alpha_{0}$}}
\put(117,11){\footnotesize{$\alpha_{1}$}}
\put(137,11){\footnotesize{$\alpha_{2}$}}
\put(151,11){\footnotesize{$\cdots$}}

\put(97,31){\footnotesize{$a$}}
\put(117,31){\footnotesize{$1$}}
\put(137,31){\footnotesize{$1$}}
\put(151,31){\footnotesize{$\cdots$}}

\put(97,51){\footnotesize{$a$}}
\put(117,51){\footnotesize{$1$}}
\put(137,51){\footnotesize{$1$}}
\put(151,51){\footnotesize{$\cdots$}}

\put(97,71){\footnotesize{$\cdots$}}
\put(117,71){\footnotesize{$\cdots$}}
\put(137,71){\footnotesize{$\cdots$}}
\put(151,71){\footnotesize{$\cdots$}}

\put(90,18){\footnotesize{$\kappa$}}
\put(90,38){\footnotesize{$1$}}
\put(90,58){\footnotesize{$1$}}
\put(91,71){\footnotesize{$\vdots$}}

\put(110,18){\footnotesize{$\frac{a\kappa}{\alpha_0}$}}
\put(110,38){\footnotesize{$1$}}
\put(110,58){\footnotesize{$1$}}
\put(111,71){\footnotesize{$\vdots$}}

\put(130,18){\footnotesize{$\frac{a\kappa}{\alpha_0\alpha_1}$}}
\put(130,38){\footnotesize{$1$}}
\put(130,58){\footnotesize{$1$}}
\put(131,71){\footnotesize{$\vdots$}}

\end{picture}
\end{center}
\caption{Weight diagram used in the Six-point Test and weight diagram of the 
$2$-variable weighted shift in Lemma \ref{prophyp}.}
\label{Figure 0}
\end{figure}

Next, given integers $i$ and $\ell $ ($\ell \geq 1$, $0\leq i\leq \ell -1$),
consider $\mathcal{H}\equiv \ell ^{2}(\mathbb{Z}_{+})=\bigvee_{j=0}^{\infty
}\{e_{j}\}.$ \ Define $\mathcal{H}_{i}:=\bigvee_{j=0}^{\infty }\{e_{\ell
j+i}\},$ so $\mathcal{H}=\bigoplus_{i=0}^{\ell -1}\mathcal{H}_{i}$.
Following the notation in \cite{CuP}, for a weight sequence $\alpha $ let 
\begin{equation}
P_{i\ell }(\alpha )\equiv \alpha (\ell :i):=\{\Pi _{m=0}^{\ell -1}\alpha
_{\ell j+i+m}\}_{j=0}^{\infty };  \label{ra}
\end{equation}%
that is, $\alpha (\ell :i)$ denotes the sequence of products of weights in
adjacent packets of size $\ell ,$ beginning with $\alpha _{i}\cdots \alpha
_{i+\ell -1}.$ \ For example, $\alpha (2:0):\alpha _{0}\alpha _{1},\alpha
_{2}\alpha _{3},\alpha _{4}\alpha _{5},\cdots $, $\alpha (2:1):\alpha
_{1}\alpha _{2},\alpha _{3}\alpha _{4},\alpha _{5}\alpha _{6},\cdots $ and $%
\alpha (3:2):\alpha _{2}\alpha _{3}\alpha _{4},\alpha _{5}\alpha _{6}\alpha
_{7},\alpha _{8}\alpha _{9}\alpha _{10},\cdots $. \ Observe that, using the
notation introduced in (\ref{sa}), $P_{i\ell }=P_{0\ell }S^{i}$. \ For a
subnormal weighted shift $W_{\alpha }$, it was proved in \cite{CuP} that $%
W_{P_{i\ell }(\alpha )}$ is also subnormal (all $\ell \geq 1$, $0\leq i\leq
\ell -1$). \ In fact, more is true.

\begin{lemma}
$($\cite{CuP}$)$\label{multiply} For $\ell \geq 1,$ and $0\leq i\leq \ell
-1, $ $W_{P_{i\ell }(\alpha )}$ is unitarily equivalent to $W_{\alpha
}^{\ell }|_{\mathcal{H}_{i}}.$ \ Therefore, $W_{\alpha }^{\ell }$ is
unitarily equivalent to $\oplus _{i=0}^{\ell -1}W_{P_{i\ell }(\alpha )}.$ \
Consequently, $W_{\alpha }^{\ell }$ is $k$-hyponormal if and only if $%
W_{P_{i\ell }(\alpha )}$ is $k$-hyponormal for each $i$ such that $0\leq
i\leq \ell -1.$ \ Moreover, if $W_{\alpha }$ is subnormal with Berger
measure $\xi _{\alpha }$, then $W_{P_{i\ell }(\alpha )}$ is subnormal with
Berger measure 
\begin{equation}
d\xi _{P_{i\ell }(\alpha )}(s)=d\xi _{P_{0\ell }S^{i}(\alpha )}(s)=\frac{%
s^{i}}{\gamma _{i}(\alpha )}d\xi _{P_{0\ell }}(s)=\frac{s^{\frac{i}{\ell }}}{%
\gamma _{i}(\alpha )}d\xi _{\alpha }(s^{\frac{1}{\ell }})\;\;(0\leq i\leq
\ell -1).  \label{rsm}
\end{equation}
\end{lemma}

\begin{example}
\label{multiply2}Let $W_{\alpha }\equiv shift(\alpha _{0},\alpha _{1},\cdots
)$ be a subnormal weighted shift, with Berger measure $\xi _{\alpha }$. \
Then $shift(\alpha _{2}\alpha _{3},\alpha _{4}\alpha _{5},...)\equiv
W_{P_{22}(\alpha )}$ is also subnormal, with Berger measure $\frac{s}{\alpha
_{0}^{2}\alpha _{1}^{2}}d\xi _{\alpha }(\sqrt{s})$. \ 
\end{example}

To produce an example of $\mathbf{T}\equiv (T_{1},T_{2})\in \mathfrak{H}_{1}$
such that $\mathbf{T}^{(2,1)}\notin \mathfrak{H}_{1}$, we start with an
example given in \cite{CLY}. \ For $0<\kappa \leq 1$, let $\alpha \equiv
\{\alpha _{n}\}_{n=0}^{\infty }$ be defined by 
\begin{equation}
\alpha _{n}:=%
\begin{cases}
\kappa \sqrt{\frac{3}{4}} & \text{if }n=0 \\ 
\frac{\sqrt{(n+1)(n+3)}}{(n+2)} & \text{if }n\geq 1.%
\end{cases}
\label{alpha}
\end{equation}%
We know that $W_{\alpha }$ is subnormal, with Berger measure 
\begin{equation*}
d\xi _{\alpha }(s):=(1-\kappa ^{2})d\delta _{0}(s)+\frac{\kappa ^{2}}{2}ds+%
\frac{\kappa ^{2}}{2}d\delta _{1}(s)\;\;\text{(\cite[Proposition 4.2]{CLY}).}
\end{equation*}%
For $0<a<1$, consider the $2$-variable weighted shift given by Figure \ref%
{Figure 0}, with $\alpha \equiv \{\alpha _{n}\}_{n=0}^{\infty }$ as above.

\begin{lemma}
(\cite{CLY})\label{prophyp} \ Let $\mathbf{T}\equiv (T_{1},T_{2})$ be the $2$%
-variable weighted shift whose weight diagram is given by Figure \ref{Figure
0}, with $0<a\leq \sqrt{\frac{1}{2}}$. \ Then\newline
(i) \ $T_{1}$ and $T_{2}$ are subnormal;\newline
(ii) \ $\mathbf{T}\in \mathfrak{H}_{1}$ if and only if $0<\kappa \leq
h_{1}(a):=\sqrt{\frac{32-48a^{4}}{{59-72a^{2}}}}$;\newline
(iii) \ $\mathbf{T}\in \mathfrak{H}_{2}$ if and only if $0<\kappa \leq
h_{2}(a):=\sqrt{\frac{81-144a^{2}}{157-360a^{2}+144a^{4}}}$;\newline
(iv) \ $\mathbf{T}\in \mathfrak{H}_{\infty }$ if and only if $0<\kappa \leq
h_{\infty }(a):=\frac{1}{\sqrt{2-a^{2}}}$.
\end{lemma}

\begin{remark}
\label{hypo}Close inspection of the proof of Lemma \ref{prophyp} reveals
that the hyponormality of the $2$-variable weighted shift $\mathbf{T}$ whose
weight diagram is given by Figure \ref{Figure 0} extends beyond the range $%
0<a\leq \sqrt{\frac{1}{2}}$. \ As a matter of fact, the hyponormality of $%
\mathbf{T}$ is controlled by the nonnegativity of the two expressions, $%
f(a):=84-95a^{2}$ and $g(a,\kappa ):=(72a^{2}-59)\kappa ^{2}+32-48a^{4}$. \
Of course, the nonnegativity of $f$ requires $a\leq \sqrt{\frac{84}{95}}$,
while to analyze the second expression we need to consider three cases: (i) $%
72a^{2}-59<0$; (ii) $72a^{2}-59=0$; and (iii) $72a^{2}-59>0$. \ In case (i), 
$g(a,\kappa )\geq 0\Leftrightarrow a^{4}\leq \frac{2}{3}$ and $\kappa
^{2}\leq \frac{32-48a^{4}}{59-72a^{2}}$; in case (ii), $a^{2}=\frac{59}{72}$
and $g(a,\kappa )=32-48(\frac{59}{72})^{2}<0$; and in case (iii), $%
g(a,\kappa )\geq 0\Leftrightarrow a^{2}>\frac{59}{72}$ and $\kappa ^{2}\geq 
\frac{32-48a^{4}}{59-72a^{2}}$. \ Now, it is easy to verify that on the
interval $(\sqrt{\frac{59}{72}},\sqrt{\frac{84}{95}}]$ the expression $\frac{%
32-48a^{4}}{59-72a^{2}}$ is always greater than $1$, and since we must have $%
\kappa \leq 1$, case (iii) cannot really happen. \ If we now observe that $%
a\leq \sqrt{\frac{84}{95}}$ is implied by the condition $a^{4}\leq \frac{2}{3%
}$, we conclude that $\mathbf{T}$ is hyponormal if and only if $a\leq \sqrt[4%
]{\frac{2}{3}}$ and $\kappa \leq \sqrt{\frac{32-48a^{4}}{59-72a^{2}}}%
=h_{1}(a)$.
\end{remark}

\begin{theorem}
\label{powhyp} Let $\mathbf{T}\equiv (T_{1},T_{2})$ be the $2$-variable
weighted shift whose weight diagram is given by Figure \ref{Figure 0}. \
Then $\mathbf{T}^{(2,1)}\equiv (T_{1}^{2},T_{2})$ is hyponormal if and only
if $0<\kappa \leq h_{21}(a):=3\sqrt{\frac{3-5a^{4}}{47-60a^{2}}}$, with $%
0<a\leq \sqrt[4]{\frac{3}{5}}$.
\end{theorem}

\textbf{Proof.} \ For $m=0,1$, let $\mathcal{H}_{m}:=\bigvee_{j=0}^{\infty
}\{e_{(2j+m,k)}:k=0,1,2,\cdots \}$. $\ $Then $\ell ^{2}(\mathbb{Z}%
_{+}^{2})\equiv \mathcal{H}_{0}\bigoplus \mathcal{H}_{1}$, and each of $%
\mathcal{H}_{0}$ and $\mathcal{H}_{1}$ reduces $T_{1}^{2}$ and $T_{2}$. \ We
can thus write 
\begin{equation*}
(T_{1}^{2},T_{2})\cong (W_{\alpha (2:0)}\oplus (I\otimes S_{a}),T_{2}|_{%
\mathcal{H}_{0}})\bigoplus (W_{\alpha (2:1)}\oplus (I\otimes U_{+}),T_{2}|_{%
\mathcal{H}_{1}}).
\end{equation*}%
By \cite[Theorem 5.2 and Remark 5.3]{CuYo1}, the second summand, $(W_{\alpha
(2:1)}\oplus (I\otimes U_{+}),T_{2}|_{\mathcal{H}_{1}})$, is subnormal. \
Thus, the hyponormality of $(T_{1}^{2},T_{2})$ is equivalent to the
hyponormality of the first summand, $(W_{\alpha (2:0)}\oplus (I\otimes
S_{a}),T_{2}|_{\mathcal{H}_{0}}).$ \ Now, to check the hyponormality of the
first summand, by Lemma \ref{joint hypo} it suffices to apply the Six-point
Test at $\mathbf{k=(}0,0)$. \ We have 
\begin{equation*}
\begin{tabular}{l}
$H_{(W_{\alpha (2:0)}\oplus (I\otimes S_{a}),T_{2}|_{\mathcal{H}_{0}})}(%
\mathbf{0})\equiv \left( 
\begin{array}{cc}
\alpha _{3}^{2}\alpha _{2}^{2}-\alpha _{1}^{2}\alpha _{0}^{2} & \frac{%
a^{2}\kappa }{\alpha _{0}\alpha _{1}}-\kappa \alpha _{0}\alpha _{1} \\ 
\frac{a^{2}\kappa }{\alpha _{0}\alpha _{1}}-\kappa \alpha _{0}\alpha _{1} & 
1-\kappa ^{2}%
\end{array}%
\right) $ \\ 
\\ 
$=\left( 
\begin{array}{cc}
\frac{9}{10}-\frac{2}{3}\kappa ^{2} & \sqrt{6}(\frac{1}{2}a^{2}-\frac{1}{3}%
\kappa ^{2}) \\ 
\sqrt{6}(\frac{1}{2}a^{2}-\frac{1}{3}\kappa ^{2}) & 1-\kappa ^{2}%
\end{array}%
\right) \geq 0$ \\ 
\\ 
$\Leftrightarrow (1-\kappa ^{2})(\frac{9}{10}-\frac{2}{3}\kappa ^{2})\geq 6(%
\frac{a^{2}}{2}-\frac{\kappa ^{2}}{3})^{2}$ \\ 
\\ 
$\Leftrightarrow \frac{9}{10}-\frac{47}{30}\kappa ^{2}-\frac{3}{2}%
a^{4}+2a^{2}\kappa ^{2}\geq 0$ \\ 
\\ 
$\Leftrightarrow h(a,\kappa ):=(60a^{2}-47)\kappa ^{2}+27-45a^{4}\geq 0$.%
\end{tabular}%
\end{equation*}%
As in Remark \ref{hypo}, three cases arise: (i) $60a^{2}-47<0$; (ii) $%
72a^{2}-59=0$; and (iii) $60a^{2}-47>0$. \ In case (i), $h(a,\kappa )\geq
0\Leftrightarrow a^{4}\leq \frac{3}{5}$ and $\kappa ^{2}\leq \frac{%
9(3-5a^{4})}{47-60a^{2}}$; in case (ii), $a^{2}=\frac{47}{60}$ and $%
h(a,\kappa )=27-45(\frac{47}{60})^{2}<0$; and in case (iii), $h(a,\kappa
)\geq 0\Leftrightarrow a^{2}>\frac{47}{60}$ and $\kappa ^{2}\geq \frac{%
27-45a^{4}}{47-60a^{2}}$. \ As before, it is easy to verify that on the
interval $(\sqrt{\frac{47}{60}},1]$ the expression $\frac{27-45a^{4}}{%
47-60a^{2}}$ is always greater than $1$, and since we must have $\kappa \leq
1$, case (iii) cannot really happen. \ We conclude that $\mathbf{T}$ is
hyponormal if and only if $a\leq \sqrt[4]{\frac{3}{5}}$ and $\kappa \leq 
\sqrt{\frac{9(3-5a^{4})}{47-60a^{2}}}\equiv h_{21}(a)$, as desired. \qed

We are now ready to formulate our first main result. \ Consider the two
functions $h_{1}$ and $h_{21}$ in Remark \ref{hypo} and Theorem \ref{powhyp}%
, respectively, restricted to the common portion of their domains, namely
the interval $(0,\sqrt[4]{\frac{3}{5}}]$. \ A calculation shows that there
exists a unique point $a_{int}\in (0,\sqrt[4]{\frac{3}{5}}]$ such that $%
h_{1}(a_{int})=h_{21}(a_{int})$; in fact, $a_{int}\cong 0.8386$. $\ $Figure %
\ref{range2} shows two regions in the $(a,\kappa )$-plane, one where $%
\mathbf{T}$ is hyponormal but $\mathbf{T}^{(2,1)}$ is not, and one where $%
\mathbf{T}^{(2,1)}$ is hyponormal but $\mathbf{T}$ is not. \ For added
emphasis, we include the graphs of $h_{2}$ and $h_{\infty }$ mentioned in
Lemma \ref{prophyp}, which are only defined on the interval $(0,\sqrt{\frac{1%
}{2}}]$. \ We thus have:

\begin{theorem}
\label{firstmain}Let $\mathbf{T}$ be the $2$-variable weighted shift whose
weight diagram is given by Figure \ref{Figure 0}. \ Then\newline
(i) \ $\mathbf{T\in }\mathfrak{H}_{1}$ and $\mathbf{T}^{(2,1)}\notin 
\mathfrak{H}_{1}\iff a_{int}<a\leq \sqrt[4]{\frac{3}{5}}$ and $%
h_{21}(a)<\kappa \leq h_{1}(a)$ (see Figure \ref{range2}).\newline
(ii) \ $\mathbf{T}\notin \mathfrak{H}_{1}$ and $\mathbf{T}^{(2,1)}\in 
\mathfrak{H}_{1}\iff 0<a<a_{int}$ and $h_{1}(a)<\kappa \leq h_{21}(a)$ (see
Figure \ref{range2}).
\end{theorem}

\setlength{\unitlength}{1mm} \psset{unit=1mm}

\begin{figure}[th]
\begin{center}
\begin{picture}(110,63)
\psline{->}(0,10)(0,58) 
\put(-1,10){$\_$} 
\put(-1,20){$\_$}
\put(-1,30){$\_$} 
\put(-1,40){$\_$} 
\put(-1,50){$\_$}
\put(-6,9){$0.6$} 
\put(-6,19){$0.7$} 
\put(-6,29){$0.8$}
\put(-6,39){$0.9$} 
\put(-6,49){$1.0$} 
\put(95,6){$a$}
\put(-4,58){$\kappa$} 

\psline(19.5,9.5)(19.5,10.5) 
\psline(39.5,9.5)(39.5,10.5) 
\psline(59.5,9.5)(59.5,10.5) 
\psline(79.5,9.5)(79.5,10.5) 

\put(-2,5){\footnotesize{$0$}} 
\put(17.3,5){\footnotesize{$0.2$}}
\put(37.3,5){\footnotesize{$0.4$}} 
\put(57.3,5){\footnotesize{$0.6$}}
\put(66.2,5){\footnotesize{$\sqrt{\frac{1}{2}}$}}
\put(77.3,5){\footnotesize{$0.8$}} 
\psline{->}(83.86,3)(83.86,9)
\put(81,1){\footnotesize{$a_{int}$}}
%\put(79.5,5){\footnotesize{$0.819$}}
\put(40,28){\footnotesize{$h_{1}$}}
\put(88,50){\footnotesize{$h_{1}$}}
\put(88,28){\footnotesize{$h_{21}$}}
\put(47,38){\footnotesize{$h_{21}$}}
\put(55,30){\footnotesize{$h_{2}$}}
\put(32,20){\footnotesize{$h_{\infty}$}}

\psline[linestyle=dashed,dash=3pt 2pt]{-}(70.7,33.3)(70.7,10)
\psline[linestyle=dashed,dash=3pt 2pt]{-}(83.86,49)(83.86,10)

\psline{->}(0,10)(99,10) 
\psline[linestyle=dashed,dash=3pt
2pt]{-}(88.1,49)(88.1,10)
\put(50,58){\footnotesize{$\mathbf{T}\in\mathfrak{H}_{1}$ and
$\mathbf{T}^{(2,1)}\notin \mathfrak{H}_{1}$ for
 $\ (a,\kappa)$ in this region}}
\psline{->}(84,56)(87,47)

\pspolygon*[linecolor=lightgray](0,23.6)(20,25.421)(40,30.50)(60,38.28)
(80,47.7)(83,49)(83,49.8)(80,49.8)(70,46)(60,41)(40,32.6)(20,27.3)(0,25.3)

\put(3,51){\footnotesize{$\mathbf{T}\notin\mathfrak{H}_{1}$ and
$\mathbf{T}^{(2,1)}\in \mathfrak{H}_{1}$ for
 $\ (a,\kappa)$ in this region}}
\psline{->}(20,49)(35,30)

%graph of h_{\infty}
\pscurve[linewidth=1pt](0,20)(20,21.428)(40,23.72)(60,28.08)(70.7,31.645)

%graph of h_{2}
\pscurve[linewidth=1pt](0,21.82)(20,22.58)(40,24.98)(60,29.565)(70.7,33.3)

%graph of h_{1}
\pscurve[linewidth=1pt](0,23.6)(20,25.421)(40,30.50)(60,38.28)(80,47.7)(83,49)(88,49.55)

%graph of h_{21}
\pscurve[linewidth=1pt](0,25.6)(20,27.7)(40,33.1)(60,41.3)(80,49.8)(83,49.8)(85,48.06)(86,45.33)(87.79,0.6)

\end{picture}
\end{center}
\caption{Graphs of $h_{1}$, $h_{21}$, $h_{2}$ and $h_{\infty }$ on the
interval $[0,\sqrt[4]{\frac{3}{5}}]$.}
\label{range2}
\end{figure}
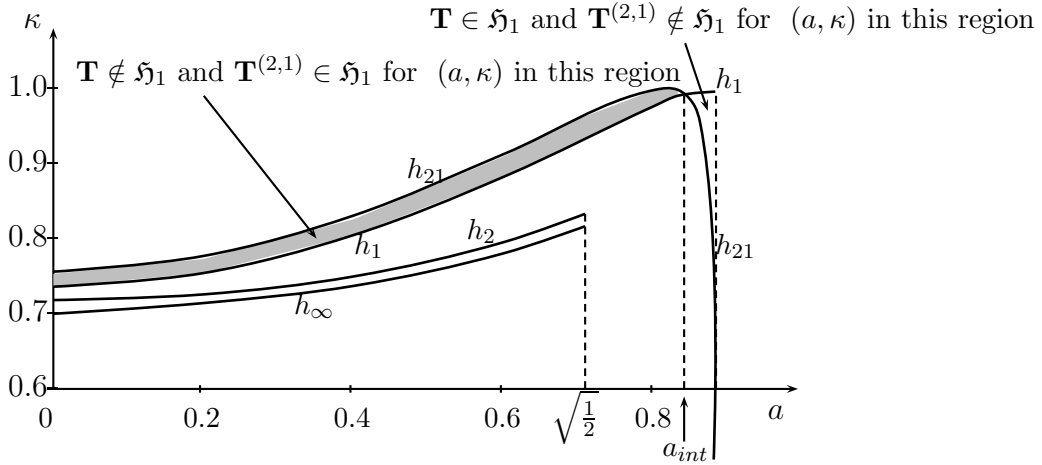

\section{A Large Class for Which $(T_{1}^{2},T_{2})\in \mathfrak{H}_{\infty
}\iff (T_{1},T_{2}^{2})\in \mathfrak{H}_{\infty }\iff (T_{1},T_{2})\in 
\mathfrak{H}_{\infty }$}

It is well known that for a single operator $T$, the subnormality of all
powers $T^{n}\;(n\geq 2)$ does not imply the hyponormality of $T$, even if $%
T $ is a unilateral weighted shift (\cite{Sta2}). \ In the multivariable
case, the analogous result is nontrivial if one further assumes that each
component is subnormal. \ To study this, we begin by recalling some useful
notation and results. \ Given a weighted shift $W_{\alpha }$, a (one-step)
backward extension of $W_{\alpha }$ is the weighted shift $W_{\alpha (x)}$,
where $\alpha (x):x,\alpha _{0},\alpha _{1},\alpha _{2},\cdots $.

\begin{lemma}
\label{backward}(Subnormal backward extension of a $1$-variable weighted
shift; cf. \cite{QHWS}, \cite[Proposition 1.5]{CuYo1}) \ Let $W_{\alpha }$
be a weighted shift whose restriction $W_{\alpha }|_{\mathcal{L}}$ to $%
\mathcal{L}:=\vee \{e_{1},e_{2},\cdots \}$ is subnormal, with Berger measure 
$\mu _{\mathcal{L}}.$ \ Then $W_{\alpha }$ is subnormal (with Berger measure 
$\mu $) if and only if\newline
(i) $\ \frac{1}{t}\in L^{1}(\mu _{\mathcal{L}})$, and\newline
(ii) $\ \alpha _{0}^{2}\leq (\left\| \frac{1}{t}\right\| _{L^{1}(\mu _{%
\mathcal{L}})})^{-1}$.\newline
In this case, $d\mu (t)=\frac{\alpha _{0}^{2}}{t}d\mu _{\mathcal{L}%
}(t)+(1-\alpha _{0}^{2}\left\| \frac{1}{t}\right\| _{L^{1}(\mu _{\mathcal{L}%
})})d\delta _{0}(t)$, where $\delta _{0}$ denotes Dirac measure at $0$. \ In
particular, $W_{\alpha }$ is never subnormal when $\mu _{\mathcal{L}%
}(\{0\})>0$. \ 
\end{lemma}

\begin{corollary}
\label{backwardcor}Let $W_{\alpha }$ be a subnormal weighted shift, let $%
\mathcal{L}_{2}:=$ $\vee \{e_{2},e_{3},\cdots \}$ and let $\mu _{\mathcal{L}%
_{2}}$ denote the Berger measure of $W_{\alpha }|_{\mathcal{L}_{2}}$. \ Then 
$\alpha _{1}$ is completely determined by $\mu _{\mathcal{L}_{2}}$, namely $%
\alpha _{1}^{2}=(\left\| \frac{1}{t}\right\| _{L^{1}(\mu _{\mathcal{L}%
_{2}})})^{-1}$. \ More generally, for $j\geq 3$ let $\mathcal{L}_{j}:=$ $%
\vee \{e_{j},e_{j+1},\cdots \}$, and let $\mu _{\mathcal{L}_{j}}$ denote the
Berger measure of $W_{\alpha }|_{\mathcal{L}_{j}}$; then $\alpha
_{j-1}=(\left\| \frac{1}{t}\right\| _{L^{1}(\mu _{\mathcal{L}_{j}})})^{-1}$.
\end{corollary}

\textbf{Proof.} \ Without loss of generality, we prove only the first
assertion. \ Since $W_{a}|_{\mathcal{L}}$ is subnormal, Lemma \ref{backward}
implies that $\alpha _{1}^{2}\leq (\left\| \frac{1}{t}\right\| _{L^{1}(\mu _{%
\mathcal{L}_{2}})})^{-1}$. \ If strict inequality occurred, then the measure 
$\mu _{\mathcal{L}}$ would have an atom at $0$, which would render the
subnormality of $W_{\alpha }$ impossible. \qed

To state the $2$-variable version of Lemma \ref{backward}, we need to recall
two notions from \cite{CuYo1}: (i) given a probability measure $\mu $ on $%
X\times Y\equiv \mathbb{R}_{+}\times \mathbb{R}_{+}$, with $\frac{1}{t}\in
L^{1}(\mu )$, the \textit{extremal measure} $\mu _{ext}$ (which is also a
probability measure) on $X\times Y$ is given by $d\mu _{ext}(s,t):=(1-\delta
_{0}(t))\frac{1}{t\left\Vert \frac{1}{t}\right\Vert _{L^{1}(\mu )}}d\mu
(s,t) $; and (ii) given a measure $\mu $ on $X\times Y$, the \textit{%
marginal measure} $\mu ^{X}$ is given by $\mu ^{X}:=\mu \circ \pi _{X}^{-1}$%
, where $\pi _{X}:X\times Y\rightarrow X$ is the canonical projection onto $%
X $. \ Thus, $\mu ^{X}(E)=\mu (E\times Y)$, for every $E\subseteq X$. \
Observe that if $\mu $ is a probability measure, then so is $\mu ^{X}$. \
For example, 
\begin{equation}
d(\xi \times \eta )_{ext}(s,t)=(1-\delta _{0}(t))\frac{1}{t\left\Vert \frac{1%
}{t}\right\Vert _{L^{1}(\eta )}}d\xi (s)d\eta (t)  \label{eq3}
\end{equation}%
and $(\xi \times \eta )^{X}=\xi $.

\begin{lemma}
\label{backext}(Subnormal backward extension of a $2$-variable weighted
shift; cf. \cite[Proposition 3.10]{CuYo1}) \ Consider the following $2$%
-variable weighted shift (see Figure \ref{compactperturbation}), and let $%
\mathcal{M}$ be the subspace of $\ell ^{2}(\mathbb{Z}_{+}^{2})$ spanned by
the canonical orthonormal basis vectors associated to indices $\mathbf{k}%
=(k_{1},k_{2})$ with $k_{1}\geq 0$ and $k_{2}\geq 1$. \ Assume that $%
\mathcal{R}_{10}(\mathbf{T})\equiv \mathbf{T}|_{\mathcal{M}}$ is subnormal
with Berger measure $\mu _{\mathcal{M}}$ and that $W_{0}:=shift(\alpha
_{00},\alpha _{10},\cdots )$ is subnormal with Berger measure $\nu $. \ Then 
$\mathbf{T}$ is subnormal if and only if \newline
(i) $\ \frac{1}{t}\in L^{1}(\mu _{\mathcal{M}})$;\newline
(ii) $\ \beta _{00}^{2}\leq (\left\| \frac{1}{t}\right\| _{L^{1}(\mu _{%
\mathcal{M}})})^{-1}$;\newline
(iii) $\ \beta _{00}^{2}\left\| \frac{1}{t}\right\| _{L^{1}(\mu _{\mathcal{M}%
})}(\mu _{\mathcal{M}})_{ext}^{X}\leq \nu $.\newline
Moreover, if $\beta _{00}^{2}\left\| \frac{1}{t}\right\| _{L^{1}(\mu _{%
\mathcal{M}})}=1,$ then $(\mu _{\mathcal{M}})_{ext}^{X}=\nu $. \ In the case
when $\mathbf{T}$ is subnormal, the Berger measure $\mu $ of $\mathbf{T}$ is
given by 
\begin{equation*}
d\mu (s,t)=\beta _{00}^{2}\left\| \frac{1}{t}\right\| _{L^{1}(\mu _{\mathcal{%
M}})}d(\mu _{\mathcal{M}})_{ext}(s,t)+(d\nu (s)-\beta _{00}^{2}\left\| \frac{%
1}{t}\right\| _{L^{1}(\mu _{\mathcal{M}})}d(\mu _{\mathcal{M}%
})_{ext}^{X}(s))d\delta _{0}(t).
\end{equation*}

\setlength{\unitlength}{1mm} \psset{unit=1mm} 
\begin{figure}[th]
\begin{center}
\begin{picture}(165,95)

\psline{->}(20,20)(85,20)
\psline(20,40)(83,40)
\psline(20,60)(83,60)
\psline(20,80)(83,80)
\psline{->}(20,20)(20,85)
\psline(40,20)(40,83)
\psline(60,20)(60,83)
\psline(80,20)(80,83)

\put(12,16){\footnotesize{$(0,0)$}}
\put(36.5,16){\footnotesize{$(1,0)$}}
\put(56.5,16){\footnotesize{$(2,0)$}}
\put(76.5,16){\footnotesize{$(3,0)$}}

\put(27,21){\footnotesize{$\alpha_{00}$}}
\put(47,21){\footnotesize{$\alpha_{10}$}}
\put(67,21){\footnotesize{$\alpha_{20}$}}
\put(81,21){\footnotesize{$\cdots$}}

\put(27,41){\footnotesize{$\alpha_{01}$}}
\put(47,41){\footnotesize{$\alpha_{11}$}}
\put(67,41){\footnotesize{$\alpha_{21}$}}
\put(81,41){\footnotesize{$\cdots$}}

\put(27,61){\footnotesize{$\alpha_{02}$}}
\put(47,61){\footnotesize{$\alpha_{12}$}}
\put(67,61){\footnotesize{$\alpha_{22}$}}
\put(81,61){\footnotesize{$\cdots$}}

\put(27,81){\footnotesize{$\cdots$}}
\put(47,81){\footnotesize{$\cdots$}}
\put(67,81){\footnotesize{$\cdots$}}

\psline{->}(50,14)(70,14)
\put(60,10){$\rm{T}_1$}
\psline{->}(10, 50)(10,70)
\put(4,60){$\rm{T}_2$}

\put(11,39){\footnotesize{$(0,1)$}}
\put(11,59){\footnotesize{$(0,2)$}}
\put(11,79){\footnotesize{$(0,3)$}}

\put(20,28){\footnotesize{$\beta_{00}$}}
\put(20,48){\footnotesize{$\beta_{01}$}}
\put(20,68){\footnotesize{$\beta_{02}$}}
\put(21,81){\footnotesize{$\vdots$}}

\put(40,28){\footnotesize{$\sqrt{\frac{\gamma_{1,1}}{\gamma_{1,0}}}$}}
\put(40,48){\footnotesize{$\sqrt{\frac{\gamma_{1,2}}{\gamma_{1,1}}}$}}
\put(40,68){\footnotesize{$\sqrt{\frac{\gamma_{1,3}}{\gamma_{1,2}}}$}}
\put(41,81){\footnotesize{$\vdots$}}

\put(60,28){\footnotesize{$\sqrt{\frac{\gamma_{2,1}}{\gamma_{2,0}}}$}}
\put(60,48){\footnotesize{$\sqrt{\frac{\gamma_{2,2}}{\gamma_{2,1}}}$}}
\put(60,68){\footnotesize{$\sqrt{\frac{\gamma_{2,3}}{\gamma_{2,2}}}$}}
\put(61,81){\footnotesize{$\vdots$}}

%\\\\\\\\\

\psline{->}(130,14)(150,14)
\put(140,10){$\rm{T}_1$}
\psline{->}(97,50)(97,70)
\put(91,60){$\rm{T}_2$}

\psline{->}(100,20)(165,20)
\psline(100,40)(163,40)
\psline(100,60)(163,60)
\psline(100,80)(163,80)

\psline{->}(100,20)(100,85)
\psline(120,20)(120,83)
\psline(140,20)(140,83)
\psline(160,20)(160,83)

\put(93,16){\footnotesize{$(0,0)$}}
\put(116,16){\footnotesize{$(1,0)$}}
\put(136,16){\footnotesize{$(2,0)$}}
\put(156,16){\footnotesize{$(3,0)$}}

\put(107,21){\footnotesize{$x_{0}$}}
\put(127,21){\footnotesize{$x_{1}$}}
\put(147,21){\footnotesize{$x_{2}$}}
\put(161,21){\footnotesize{$\cdots$}}

\put(107,41){\footnotesize{$x$}}
\put(127,41){\footnotesize{$\alpha_{1}$}}
\put(147,41){\footnotesize{$\alpha_{2}$}}
\put(161,41){\footnotesize{$\cdots$}}

\put(107,63){\footnotesize{$\frac{x\beta_{1}}{y_{1}}$}}
\put(127,61){\footnotesize{$\alpha_{1}$}}
\put(147,61){\footnotesize{$\alpha_{2}$}}
\put(161,61){\footnotesize{$\cdots$}}

\put(107,81){\footnotesize{$\cdots$}}
\put(127,81){\footnotesize{$\cdots$}}
\put(147,81){\footnotesize{$\cdots$}}
\put(161,81){\footnotesize{$\cdots$}}

\put(100,28){\footnotesize{$y_{0}$}}
\put(100,48){\footnotesize{$y_{1}$}}
\put(100,68){\footnotesize{$y_{2}$}}
\put(101,81){\footnotesize{$\vdots$}}

\put(120,28){\footnotesize{$\frac{xy_{0}}{x_{0}}$}}
\put(120,48){\footnotesize{$\beta_{1}$}}
\put(120,68){\footnotesize{$\beta_{2}$}}
\put(121,81){\footnotesize{$\vdots$}}

\put(140,28){\footnotesize{$\frac{xy_{0}\alpha_{1}}{x_{0}x_{1}}$}}
\put(140,48){\footnotesize{$\beta_{1}$}}
\put(140,68){\footnotesize{$\beta_{2}$}}
\put(141,81){\footnotesize{$\vdots$}}

\end{picture}
\end{center}
\caption{Weight diagram of the $2$-variable weighted shift in Lemma \ref%
{backext} and weight diagram of a $2$-variable weighted shift with $\mathcal{%
R}_{11}(\mathbf{T})\protect\cong(I\otimes W_{\protect\alpha } ,W_{\protect%
\beta }\otimes I)$, respectively.}
\label{compactperturbation}
\end{figure}
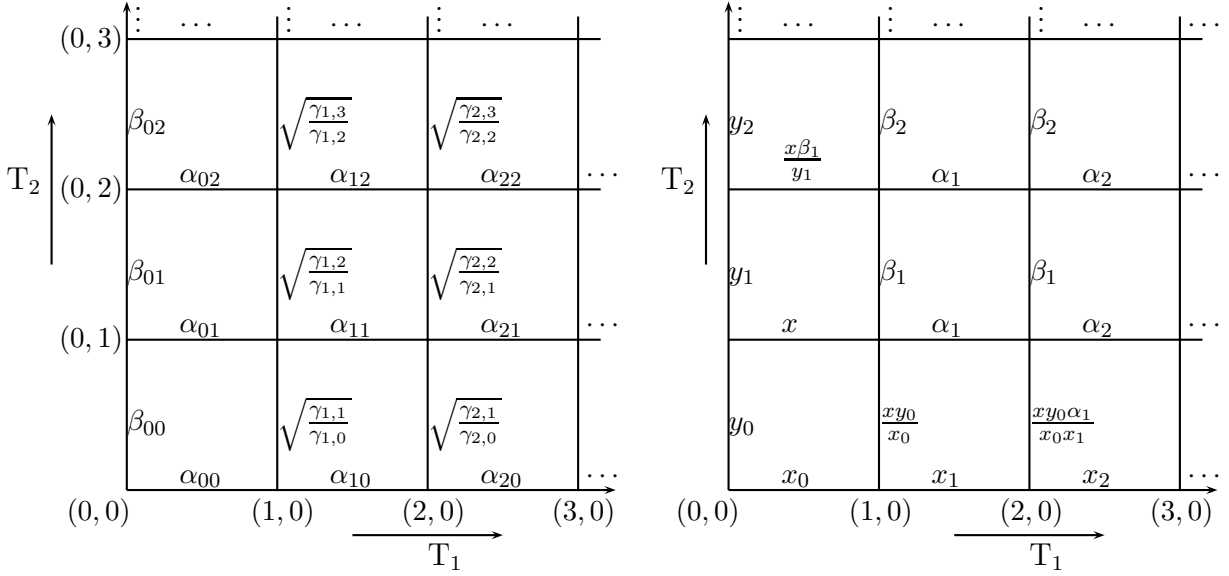
\end{lemma}

\begin{definition}
\label{tc} (i) \ The \textit{core} of a $2$-variable weighted shift $\mathbf{%
T}$ is the restriction of $\mathbf{T}$ to $\mathcal{M}_{1}\cap \mathcal{N}%
_{1}$, in symbols, $c(\mathbf{T}):=\mathbf{T}|_{\mathcal{M}_{1}\cap \mathcal{%
N}_{1}}\equiv \mathcal{R}_{11}(\mathbf{T})$.\newline
(ii) \ A $2$-variable weighted shift $\mathbf{T}$ is said to be of \textit{%
tensor form} if $\mathbf{T}\cong (I\otimes W_{\alpha },W_{\beta }\otimes I)$%
. $\ $When $\mathbf{T}$ is subnormal, this is equivalent to requiring that
the Berger measure be a Cartesian product $\xi \times \eta $.\newline
(iii) $\ $The class of all $2$-variable weighted shifts $\mathbf{T}\in 
\mathfrak{H}_{0}$ whose cores are of tensor form will be denoted by $%
\mathcal{TC}$, that is, $\mathcal{TC}:=\{\mathbf{T}\in \mathfrak{H}_{0}:c(%
\mathbf{T})\text{ is of tensor form}\}$.\newline
(iv) \ For each $\mathbf{k}\equiv (k_{1},k_{2})\in \mathbb{Z}_{+}^{2}$, let $%
A_{\mathbf{k}}:=\{\mathbf{T}\in \mathfrak{H}_{0}:\mathcal{R}_{k_{1}k_{2}}(%
\mathbf{T})\in \mathcal{TC}\}$.
\end{definition}

Observe that for $\mathbf{k},\mathbf{m}\in \mathbb{Z}_{+}^{2}$ with $\mathbf{%
k}\leq \mathbf{m}$ (i.e., $\mathbf{m}-\mathbf{k}\in \mathbb{Z}_{+}^{2}$), we
have $A_{\mathbf{k}}\subseteq A_{\mathbf{m}}.$ \ Thus, the collection $\{A_{%
\mathbf{k}}\}$ forms an ascending chain with respect to set inclusion and
the partial order induced by $\mathbb{Z}_{+}^{2}$. \ Moreover $\mathcal{TC}%
=A_{\mathbf{0}}\subseteq A_{\mathbf{k}}$ for all $\mathbf{k}\in \mathbb{Z}%
_{+}^{2}.$ \ All $2$-variable weighted shifts considered in \cite{CuYo1}, %
\cite{CuYo2}, \cite{CuYo3} and \cite{CLY} are in $\mathcal{TC}$. \ Thus, $%
\mathcal{TC}$ is a rather large class; as a matter of fact, much more is
true. \ The following theorem shows that an outer propagation phenomena
occurs for $\mathcal{TC}$.

\begin{theorem}
For all $\mathbf{k}\in \mathbb{Z}_{+}^{2},$ $A_{\mathbf{k}}=\mathcal{TC}$.
\end{theorem}

\textbf{Proof.} \ Since we always have $\mathcal{TC}\subseteq A_{\mathbf{k}}$%
, we prove the reverse inclusion. \ Without loss of generality, it is enough
to show that if $\mathbf{T}\in \mathfrak{H}_{0}$ and $\mathbf{T}|_{\mathcal{M%
}_{2}\cap \mathcal{N}_{2}}$ is of tensor form, then $c(\mathbf{T})$ is of
tensor form. \ If $\mathbf{T}|_{\mathcal{M}_{2}\cap \mathcal{N}_{2}}$ is of
tensor form, then $shift(\beta _{22},\beta _{23},\cdots )=shift(\beta
_{k_{1}2},\beta _{k_{1}3},\cdots )$ for all $k_{1}\geq 2$. \ The
subnormality of $T_{2}$ then implies that $shift(\beta _{k_{1}0},\beta
_{k_{1}1},\cdots )$ is subnormal for all $k_{1}\geq 2$. \ By Corollary \ref%
{backwardcor}, we have $\beta _{k_{1}1}=\sqrt{(\left\| \frac{1}{t}\right\|
_{L^{1}(\xi _{k_{1}})})^{-1}}$ $(k_{1}\geq 2)$, where $\xi _{k_{1}}$ is the
Berger measure of $shift(\beta _{k_{1}2},\beta _{k_{1}3},\cdots )$. \ Thus, $%
shift(\beta _{21},\beta _{22},\cdots )=shift(\beta _{k_{1}1},\beta
_{k_{1}2},\cdots )$ for all $k_{1}\geq 2$. \ Now, since $\beta _{21}=\beta
_{k_{1}1}$ (all $k_{1}\geq 2$), the commutativity of $T_{1}$ and $T_{2}$
implies $\alpha _{k_{1}2}=\alpha _{k_{1}1}$ for all $k_{1}\geq 2$. \ Thus, $%
shift(\alpha _{21},\alpha _{31},\cdots )=shift(\alpha _{2k_{2}},\alpha
_{3k_{2}},\cdots )$ for all $k_{2}\geq 1$. \ By the subnormality of $T_{1}$
and Lemma \ref{backward}, we have $shift(\alpha _{11},\alpha _{21},\cdots
)=shift(\alpha _{1k_{2}},\alpha _{2k_{2}},\cdots )$ for all $k_{2}\geq 1$. \
Therefore $c(\mathbf{T})$ is of tensor form. \qed

We now consider the $2$-variable weighted shift given by Figure \ref%
{compactperturbation}, where $W_{x}:=shift(x_{0},x_{1},\cdots )$ and $%
W_{y}:=shift(y_{0},y_{1},\cdots )$ are subnormal with Berger measures $\mu
_{y}$ and $\mu _{x}$, respectively. \ Further, we let $W_{\alpha
}:=shift(\alpha _{1},\alpha _{2},\cdots )$ and $W_{\beta }:=shift(\beta
_{1},\beta _{2},\cdots )$ be subnormal with Berger measures $\xi $ and $\eta 
$, respectively, and we let $r:=\left\| \frac{1}{s}\right\| _{L^{1}(\xi
)}\in (0,\infty ]$ and $d\widetilde{\xi }(s):=\frac{d\xi (s)}{s}$. $\ $We
then have:

\begin{theorem}
\label{thm0}Let $\mathbf{T}\equiv \mathbf{(}T_{1},T_{2})\in \mathcal{TC}.$ \
Then $\mathcal{R}_{10}(\mathbf{T})\in \mathfrak{H}_{\infty }$ if and only if 
$x^{2}r\eta \leq \left( \mu _{y}\right) _{1}$. \ In this case, the Berger
measure of $\mathcal{R}_{10}(\mathbf{T})$ is $x^{2}\widetilde{\xi }\times
\eta +\delta _{0}\times (\left( \eta _{y}\right) _{1}-x^{2}r\eta )\text{,}$
where $\left( \eta _{y}\right) _{1}$ is the Berger measure of the subnormal
shift $shift(y_{1},y_{2},\cdots )$.
\end{theorem}

\textbf{Proof.} \ This is a straightforward application of Lemma \ref%
{backext}, if we think of $\mathcal{R}_{10}(\mathbf{T})$ as the backward
extension of $c(\mathbf{T)}$ (in the $s$ direction). \qed

\begin{proposition}
Let $\mathbf{T}\equiv \mathbf{(}T_{1},T_{2})\in \mathcal{TC}.$ \ Then $%
\left\| \frac{1}{t}\right\| _{L^{1}(\left( \eta _{y}\right)
_{1})}=y_{1}^{2}\left\| \frac{1}{t}\right\| _{L^{1}(\left( \eta _{y}\right)
_{2}^{2})}$, where $\left( \eta _{y}\right) _{1}$ (resp. $\left( \eta
_{y}\right) _{2}^{2}$) is the Berger measure of $shift(y_{1},y_{2},\cdots )$
(resp. $shift(y_{2}y_{3},y_{4}y_{5},\cdots )$). \ Moreover, $\left\| \frac{1%
}{t}\right\| _{L^{1}(\eta )}=\beta _{1}^{2}\left\| \frac{1}{t}\right\|
_{L^{1}(\eta _{1}^{2})}$, where $\eta _{1}^{2}$ is the Berger measure of $%
shift(\beta _{2}\beta _{3},\beta _{4}\beta _{5},\cdots )$.
\end{proposition}

\textbf{Proof.} \ Since $shift(y_{0},y_{1},\cdots )$ has Berger measure $%
\eta _{y}$, we have $\left( d\eta _{y}\right) _{1}=\frac{t}{y_{0}^{2}}d\eta
_{y}(t)$; moreover, the Berger measure of $shift(y_{2},y_{3},\cdots )$ is $%
\left( d\eta _{y}\right) _{2}(t)=\frac{t^{2}}{y_{0}^{2}y_{1}^{2}}d\eta
_{y}(t)$. \ Thus by Lemma \ref{multiply}, \newline
$shift(y_{2}y_{3},y_{4}y_{5},\cdots )$ has Berger measure $\left( d\eta
_{y}\right) _{2}^{2}\equiv \frac{t}{y_{0}^{2}y_{1}^{2}}d\eta _{y}(\sqrt{t})$%
. \ Observe that 
\begin{equation*}
\left\| \frac{1}{t}\right\| _{L^{1}(\left( \eta _{y}\right) _{1})}=\frac{1}{%
y_{0}^{2}}=\int_{0}^{A}\frac{1}{y_{0}^{2}}d\eta _{y}(t)=\frac{1}{y_{0}^{2}}%
\int_{0}^{A^{2}}d\eta _{y}(\sqrt{t})=\frac{1}{y_{0}^{2}}\int_{0}^{A^{2}}%
\frac{y_{0}^{2}y_{1}^{2}}{t}d\left( \eta _{y}\right)
_{2}^{2}=y_{1}^{2}\left\| \frac{1}{t}\right\| _{L^{1}(\left( \eta
_{y}\right) _{2}^{2})},
\end{equation*}%
where $A:=\left\| shift(y_{0},y_{1},\cdots )\right\| ^{2}$. \ Thus, we get $%
\left\| \frac{1}{t}\right\| _{L^{1}(\left( \eta _{y}\right)
_{1})}=y_{1}^{2}\left\| \frac{1}{t}\right\| _{L^{1}(\left( \eta _{y}\right)
_{2}^{2})}$, as desired.

\ Next, we observe that $d\eta _{1}(t)\equiv \frac{t}{\beta _{1}^{2}}d\eta
(t)$ is the Berger measure of $shift(\beta _{2},\beta _{3},\cdots )$ and $%
d\eta _{1}^{2}(t)\equiv \frac{\sqrt{t}}{\beta _{1}^{2}}d\eta (\sqrt{t})$ is
the Berger measure of $shift(\beta _{2}\beta _{3},\beta _{4}\beta
_{5},\cdots )$. \ Let $B:=\left\| shift(\beta _{0},\beta _{1},\cdots
)\right\| ^{2}$; we then have 
\begin{equation*}
\left\| \frac{1}{t}\right\| _{L^{1}(\eta )}=\int_{0}^{B}\frac{1}{t}d\eta
(t)=\int_{0}^{B^{2}}\frac{1}{\sqrt{t}}d\eta (\sqrt{t})=\beta
_{1}^{2}\int_{0}^{B^{2}}\frac{1}{t}\frac{\sqrt{t}}{\beta _{1}^{2}}d\eta (%
\sqrt{t})=\beta _{1}^{2}\left\| \frac{1}{t}\right\| _{L^{1}(\eta _{1}^{2})},
\end{equation*}%
\ as desired. \qed

We next recall that $\mathbf{(}T_{1},T_{2}^{2})$ can be regarded as the
orthogonal direct sum of two $2$-variable weighted shifts. \ For $m=0,1$, let%
\begin{equation*}
\mathcal{H}^{m}:=\bigvee_{k=0}^{\infty }\{e_{(j,2k+m)}:j=0,1,2,\cdots \}.
\end{equation*}%
Then $\ell ^{2}(\mathbb{Z}_{+}^{2})\equiv \mathcal{H}^{0}\oplus \mathcal{H}%
^{1}$ and each of $\mathcal{H}^{0}$ and $\mathcal{H}^{1}$ reduces $T_{1}$
and $T_{2}$. \ Thus, $\mathbf{(}T_{1},T_{2}^{2})$ is subnormal if and only
if each of $\mathbf{(}T_{1},T_{2}^{2})|_{\mathcal{H}^{0}}$ and $\mathbf{(}%
T_{1},T_{2}^{2})|_{\mathcal{H}^{1}}$ is subnormal. \ The weight diagrams of
these $2$-variable weighted shifts are shown in Figure \ref{double}. %
\setlength{\unitlength}{1mm} \psset{unit=1mm}

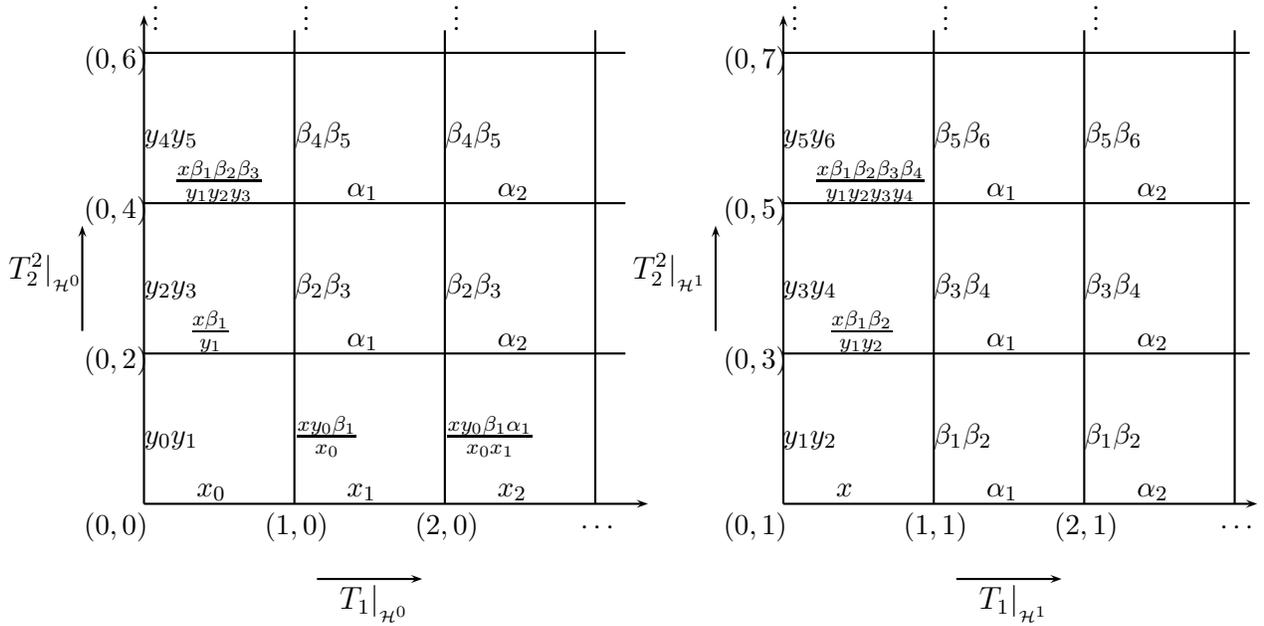
\begin{figure}[th]
\begin{center}
\begin{picture}(160,90)

\psline{->}(15,20)(82,20)
\psline(15,40)(79,40)
\psline(15,60)(79,60)
\psline(15,80)(79,80)
\psline{->}(15,20)(15,85)
\psline(35,20)(35,83)
\psline(55,20)(55,83)
\psline(75,20)(75,83)

\put(7,16){\footnotesize{$(0,0)$}}
\put(31,16){\footnotesize{$(1,0)$}}
\put(51,16){\footnotesize{$(2,0)$}}
\put(73,16){\footnotesize{$\cdots$}}

\put(22,21){\footnotesize{$x_{0}$}}
\put(42,21){\footnotesize{$x_{1}$}}
\put(62,21){\footnotesize{$x_{2}$}}

\put(21,42){\footnotesize{$\frac{x\beta_{1}}{y_{1}}$}}
\put(42,41){\footnotesize{$\alpha_{1}$}}
\put(62,41){\footnotesize{$\alpha_{2}$}}

\put(19,62){\footnotesize{$\frac{x\beta_{1}\beta_{2}\beta_{3}}{y_{1}y_{2}y_{3}}$}}
\put(42,61){\footnotesize{$\alpha_{1}$}}
\put(62,61){\footnotesize{$\alpha_{2}$}}

\psline{->}(38,10)(52,10)
\put(41,6){$T_1|_{_{\mathcal{H}^{0}}}$}

\put(7,38){\footnotesize{$(0,2)$}}
\put(7,58){\footnotesize{$(0,4)$}}
\put(7,78){\footnotesize{$(0,6)$}}

\psline{->}(6.8,43)(6.8,57)
\put(-3,50){$T_{2}^2|_{_{\mathcal{H}^{0}}}$}

\put(15,28){\footnotesize{$y_{0}y_{1}$}}
\put(15,48){\footnotesize{$y_{2}y_{3}$}}
\put(15,68){\footnotesize{$y_{4}y_{5}$}}
\put(16,83){\footnotesize{$\vdots$}}

\put(35,28){\footnotesize{$\frac{xy_{0}\beta_{1}}{x_{0}}$}}
\put(35,48){\footnotesize{$\beta_{2}\beta_{3}$}}
\put(35,68){\footnotesize{$\beta_{4}\beta_{5}$}}
\put(36,83){\footnotesize{$\vdots$}}

\put(55,28){\footnotesize{$\frac{xy_{0}\beta_{1}\alpha_{1}}{x_{0}x_{1}}$}}
\put(55,48){\footnotesize{$\beta_{2}\beta_{3}$}}
\put(55,68){\footnotesize{$\beta_{4}\beta_{5}$}}
\put(56,83){\footnotesize{$\vdots$}}

%\\\\\

\psline{->}(100,20)(163,20)
\psline(100,40)(162,40)
\psline(100,60)(162,60)
\psline(100,80)(162,80)
\psline{->}(100,20)(100,85)
\psline(120,20)(120,83)
\psline(140,20)(140,83)
\psline(160,20)(160,83)

\put(92,16){\footnotesize{$(0,1)$}}
\put(116,16){\footnotesize{$(1,1)$}}
\put(136,16){\footnotesize{$(2,1)$}}
\put(158,16){\footnotesize{$\cdots$}}

\put(107,21){\footnotesize{$x$}}
\put(127,21){\footnotesize{$\alpha_{1}$}}
\put(147,21){\footnotesize{$\alpha_{2}$}}

\put(106,42){\footnotesize{$\frac{x\beta_{1}\beta_{2}}{y_{1}y_{2}}$}}
\put(127,41){\footnotesize{$\alpha_{1}$}}
\put(147,41){\footnotesize{$\alpha_{2}$}}

\put(104,62){\footnotesize{$\frac{x\beta_{1}\beta_{2}\beta_{3}\beta_{4}}{y_{1}y_{2}y_{3}y_{4}}$}}
\put(127,61){\footnotesize{$\alpha_{1}$}}
\put(147,61){\footnotesize{$\alpha_{2}$}}

\psline{->}(123,10)(137,10) \put(126,6){$T_1|_{_{\mathcal{H}^{1}}}$}

\put(92,38){\footnotesize{$(0,3)$}}
\put(92,58){\footnotesize{$(0,5)$}}
\put(92,78){\footnotesize{$(0,7)$}}

\psline{->}(91,43)(91,57)
\put(80,50){$T_{2}^2|_{_{\mathcal{H}^{1}}}$}

\put(100,28){\footnotesize{$y_{1}y_{2}$}}
\put(100,48){\footnotesize{$y_{3}y_{4}$}}
\put(100,68){\footnotesize{$y_{5}y_{6}$}}
\put(101,83){\footnotesize{$\vdots$}}

\put(120,28){\footnotesize{$\beta_{1}\beta_{2}$}}
\put(120,48){\footnotesize{$\beta_{3}\beta_{4}$}}
\put(120,68){\footnotesize{$\beta_{5}\beta_{6}$}}
\put(121,83){\footnotesize{$\vdots$}}

\put(140,28){\footnotesize{$\beta_{1}\beta_{2}$}}
\put(140,48){\footnotesize{$\beta_{3}\beta_{4}$}}
\put(140,68){\footnotesize{$\beta_{5}\beta_{6}$}}
\put(141,83){\footnotesize{$\vdots$}}

\end{picture}
\end{center}
\caption{Weight diagrams of $\mathbf{(}T_{1},T_{2}^{2})|_{_{\mathcal{H}%
^{0}}} $ and $\mathbf{(}T_{1},T_{2}^{2})|_{_{\mathcal{H}^{1}}}$ in the proof
of Proposition \ref{pro1} and Theorem \ref{Thm1}.}
\label{double}
\end{figure}

We first focus on $\mathbf{(}T_{1},T_{2}^{2})|_{\mathcal{H}^{1}}$:

\begin{proposition}
\label{pro1}Let $\mathbf{T}\equiv \mathbf{(}T_{1},T_{2})\in \mathcal{TC}$. \
Then $(T_{1},T_{2}^{2})|_{\mathcal{H}_{1}}$ is subnormal if and only if $%
\mathcal{R}_{10}(\mathbf{T})$ is subnormal.
\end{proposition}

\textbf{Proof.} \ First, recall that $shift(y_{0},y_{1},y_{2},\cdots )$ has
Berger measure $\eta _{y}$, that $d\left( \eta _{y}\right) _{1}(t)=\frac{t}{%
y_{0}^{2}}d\eta _{y}(t)$ and that $d\left( \eta _{y}\right) _{2}(t)=\frac{%
t^{2}}{y_{0}^{2}y_{1}^{2}}d\eta _{y}(t)$. \ Now, Theorem \ref{thm0} states
that%
\begin{equation*}
(T_{1},T_{2})|_{\mathcal{M}_{1}}\text{ is subnormal }\Leftrightarrow
x^{2}r\eta \leq \left( \eta _{y}\right) _{1}\text{.}
\end{equation*}%
On the other hand, Theorem \ref{thm0} (applied to $(T_{1},T_{2}^{2})|_{%
\mathcal{H}_{1}}$) says that%
\begin{equation*}
(T_{1},T_{2}^{2})|_{\mathcal{H}_{1}}\text{ is subnormal }\iff x^{2}r\eta
^{2}\leq \left( \eta _{y}\right) _{1}^{2},
\end{equation*}%
and if $\mathbf{(}T_{1},T_{2})|_{\mathcal{H}_{1}}$ is subnormal, its Berger
measure is $x^{2}\widetilde{\xi }\times \eta ^{2}+\delta _{0}\times (\left(
\eta _{y}\right) _{1}^{2}-x^{2}r\eta ^{2})\text{,}$ where $\left( \eta
_{y}\right) _{1}^{2}$ is the Berger measure of $shift(y_{1}y_{2},y_{3}y_{4},%
\cdots )$ and $\eta ^{2}$ is the Berger measure of $W_{\beta }:=shift(\beta
_{1}\beta _{2},\beta _{3}\beta _{4},\cdots )$. \ By observing that 
\begin{equation*}
x^{2}r\eta ^{2}\leq \left( \eta _{y}\right) _{1}^{2}\Leftrightarrow
x^{2}rd\eta (\sqrt{t})\leq d\left( \eta _{y}\right) _{1}(\sqrt{t}%
)\Leftrightarrow x^{2}rd\eta (t)\leq d\left( \eta _{y}\right) _{1}(t),
\end{equation*}%
we obtain the desired result. \qed

\begin{theorem}
\label{Thm1}Let $\mathbf{T}\equiv (T_{1},T_{2})\in \mathcal{TC}.$ Then 
\begin{equation*}
(T_{1},T_{2}^{2})\in \mathfrak{H}_{\infty }\iff (T_{1}^{2},T_{2})\in 
\mathfrak{H}_{\infty }\iff (T_{1},T_{2})\in \mathfrak{H}_{\infty }.
\end{equation*}
\end{theorem}

\begin{corollary}
\label{cor1}Let $\mathbf{T}\equiv (T_{1},T_{2})\in \mathcal{TC}.$ If $%
(T_{1},T_{2}^{2}),\ (T_{1}^{2},T_{2})\in \mathfrak{H}_{\infty }$, then $%
(T_{1},T_{2})\in \mathfrak{H}_{\infty }$.
\end{corollary}

In view of Corollary \ref{cor1}, the following conjecture for $2$-variable
weighted shifts seems natural.

\begin{conjecture}
\label{conjecture} If $(T_{1},T_{2}^{2}),(T_{1}^{2},T_{2})\in \mathfrak{H}%
_{\infty }$, then $(T_{1},T_{2})\in \mathfrak{H}_{\infty }$.
\end{conjecture}

\textbf{Proof of Theorem \ref{Thm1}.} \ Clearly, it is enough to show that $%
(T_{1},T_{2}^{2})\in \mathfrak{H}_{\infty }\Rightarrow (T_{1},T_{2})\in 
\mathfrak{H}_{\infty }.$ \ Since $(T_{1},T_{2}^{2})\in \mathfrak{H}_{\infty
}\Rightarrow \mathbf{(}T_{1},T_{2}^{2})|_{\mathcal{H}^{0}}\in \mathfrak{H}%
_{\infty }$, our strategy consists of first characterizing the subnormality
of $\mathbf{T}$ and of $\mathbf{(}T_{1},T_{2}^{2})|_{\mathcal{H}^{0}}$ in
terms of the given parameters ($y_{0}$, $\nu $, etc), and then establishing
the desired implication at the parameter level. \ That is, we will show that 
$\mathbf{(}T_{1},T_{2}^{2})|_{\mathcal{H}^{0}}\in \mathfrak{H}_{\infty
}\Rightarrow \mathbf{T}\in \mathfrak{H}_{\infty }$ using their parametric
characterizations. \ Proposition \ref{pro1} will help us characterize the
subnormality of $\mathbf{T}$. \ Recall that $(T_{1},T_{2}^{2})|_{\mathcal{H}%
_{1}}$ is subnormal if and only if $(T_{1},T_{2})|_{\mathcal{M}_{1}}$ is
subnormal, and in that case the Berger measure of $\mathbf{(}T_{1},T_{2})|_{%
\mathcal{M}_{1}}$ is%
\begin{equation*}
\mu _{\mathcal{M}}=x^{2}\widetilde{\xi }\times \eta +\delta _{0}\times
(\left( \eta _{y}\right) _{1}-x^{2}r\eta )\text{.}
\end{equation*}%
We then have 
\begin{equation}
\left\| \frac{1}{t}\right\| _{L^{1}(\mu _{\mathcal{M}})}=\int \frac{1}{t}%
d\mu _{\mathcal{M}}(s,t)=x^{2}r\left\| \frac{1}{t}\right\| _{L^{1}(\eta
)}+\int \frac{1}{t}d\left( \eta _{y}\right) _{1}(t)-x^{2}r\left\| \frac{1}{t}%
\right\| _{L^{1}(\eta )}=\left\| \frac{1}{t}\right\| _{L^{1}(\left( \eta
_{y}\right) _{1})}.  \label{eq10}
\end{equation}%
Thus, we get 
\begin{eqnarray*}
d(\mu _{\mathcal{M}})_{ext}(s,t) &=&d\{x^{2}\widetilde{\xi }\times \eta
+\delta _{0}\times (\left( \eta _{y}\right) _{1}-x^{2}r\eta )\}_{ext}(s,t) \\
&& \\
&=&\frac{1}{t\left\| \frac{1}{t}\right\| _{L^{1}(\mu _{\mathcal{M}})}}%
\{x^{2}d\widetilde{\xi }(s)d\eta (t)+d\delta _{0}(s)(d\left( \eta
_{y}\right) _{1}(t)-x^{2}rd\eta (t))\} \\
&& \\
&=&\frac{1}{\left\| \frac{1}{t}\right\| _{L^{1}(\mu _{\mathcal{M}})}}\left\{
x^{2}d\widetilde{\xi }(s)\frac{d\eta (t)}{t}+d\delta _{0}(s)\left( \frac{%
d\left( \eta _{y}\right) _{1}(t)}{t}-x^{2}r\frac{d\eta (t)}{t}\right)
\right\} .
\end{eqnarray*}%
From (\ref{eq10}), it follows that 
\begin{equation}
(\mu _{\mathcal{M}})_{ext}^{X}=\left( \frac{x^{2}\left\| \frac{1}{t}\right\|
_{L^{1}(\eta )}}{\left\| \frac{1}{t}\right\| _{L^{1}(\left( \eta _{y}\right)
_{1})}}\right) \widetilde{\xi }+\left( 1-\frac{x^{2}r\left\| \frac{1}{t}%
\right\| _{L^{1}(\eta )}}{\left\| \frac{1}{t}\right\| _{L^{1}(\left( \eta
_{y}\right) _{1})}}\right) \delta _{0}.  \label{eq10p}
\end{equation}%
If we let $\varphi $ denote the right-hand side in (\ref{eq10p}), it follows
that 
\begin{eqnarray}
(T_{1},T_{2})\text{ is subnormal } &\Leftrightarrow &y_{0}^{2}\left\| \frac{1%
}{t}\right\| _{L^{1}(\mu _{\mathcal{M}})}(\mu _{\mathcal{M}})_{ext}^{X}\leq
\nu \;\;\text{(by Lemma \ref{backext})}  \label{eq11} \\
&&  \notag \\
&\Leftrightarrow &y_{0}^{2}\left\| \frac{1}{t}\right\| _{L^{1}(\left( \eta
_{y}\right) _{1})}\varphi \leq \mu _{x}\;\;\text{(using (\ref{eq10}).} 
\notag
\end{eqnarray}%
We have thus characterized the subnormality of $\mathbf{T}$. \ \newline
We now consider the $2$-variable weighted shift $(T_{1},T_{2}^{2})|_{%
\mathcal{H}^{0}}$ and the associated subspace $\mathcal{HM}:=\mathbf{\vee }%
\{e_{\mathbf{k}}\in \mathcal{H}^{0}:k_{2}\geq 1\}$. \ Observe that $%
(T_{1},T_{2}^{2})|_{\mathcal{H}^{0}}$ can be regarded as a backward
extension of $(T_{1},T_{2}^{2})|_{\mathcal{HM}}$, and that the latter is
subnormal with Berger measure%
\begin{equation*}
\theta :=\frac{x^{2}\beta _{1}^{2}}{y_{1}^{2}}\widetilde{\xi }\times \eta
_{1}^{2}+\delta _{0}\times \left( \left( \eta _{y}\right) _{2}^{2}-\frac{%
x^{2}r\beta _{1}^{2}}{y_{1}^{2}}\eta _{1}^{2}\right) ,
\end{equation*}%
where $\eta _{1}^{2}$ (resp. $\left( \eta _{y}\right) _{2}^{2}$) is the
Berger measure of $shift(\beta _{2}\beta _{3},\beta _{4}\beta _{5},\cdots )$
(resp. $shift(y_{2}y_{3},y_{4}y_{5},\cdots ))$. We then have 
\begin{equation}
\left\| \frac{1}{t}\right\| _{L^{1}(\mu _{\mathcal{HM}})}=\left\| \frac{1}{t}%
\right\| _{L^{1}(\left( \eta _{y}\right) _{2}^{2})},  \label{eq12}
\end{equation}%
and%
\begin{eqnarray*}
d(\mu _{\mathcal{HM}})_{ext}(s,t) &=&d\theta _{ext}(s,t) \\
&& \\
&=&\frac{1}{t\left\| \frac{1}{t}\right\| _{L^{1}(\mu _{\mathcal{HM}})}}%
d\theta (s,t).
\end{eqnarray*}%
From (\ref{eq12}), we have%
\begin{equation}
(\mu _{\mathcal{HM}})_{ext}^{X}=\frac{1}{\left\| \frac{1}{t}\right\|
_{L^{1}(\left( \eta _{y}\right) _{2}^{2})}}\left\{ \frac{x^{2}\beta _{1}^{2}%
}{y_{1}^{2}}\left\| \frac{1}{t}\right\| _{L^{1}(\eta _{1}^{2})}\widetilde{%
\xi }+(\left\| \frac{1}{t}\right\| _{L^{1}(\left( \eta _{y}\right)
_{2}^{2})}-\frac{x^{2}r\beta _{1}^{2}}{y_{1}^{2}}\left\| \frac{1}{t}\right\|
_{L^{1}(\eta _{1}^{2})})\delta _{0}\right\} .  \label{eq12p}
\end{equation}%
If we now let $\psi $ denote the expression in braces in the right-hand side
of (\ref{eq12p}), Lemma \ref{backext} combined with (\ref{eq12}) imply that 
\begin{equation}
\begin{tabular}{l}
$(T_{1},T_{2}^{2})|_{\mathcal{H}^{0}}\text{ is subnormal }$ \\ 
\\ 
$\Leftrightarrow y_{0}^{2}y_{1}^{2}\left\| \frac{1}{t}\right\| _{L^{1}(\mu _{%
\mathcal{HM}})}(\mu _{\mathcal{HM}})_{ext}^{X}\leq \nu \Leftrightarrow
y_{0}^{2}y_{1}^{2}\psi \leq \mu _{x}$ \\ 
\\ 
$\Leftrightarrow y_{0}^{2}\left\{ x^{2}\beta _{1}^{2}\left\| \frac{1}{t}%
\right\| _{L^{1}(\eta _{1}^{2})}\widetilde{\xi }+\left( y_{1}^{2}\left\| 
\frac{1}{t}\right\| _{L^{1}(\left( \eta _{y}\right) _{2}^{2})}-x^{2}r\beta
_{1}^{2}\left\| \frac{1}{t}\right\| _{L^{1}(\eta _{1}^{2})}\right) \delta
_{0}\right\} \leq \mu _{x}$.%
\end{tabular}
\label{eq13}
\end{equation}%
Observe that%
\begin{equation}
\begin{tabular}{l}
$y_{0}^{2}\left\| \frac{1}{t}\right\| _{L^{1}(\left( \eta _{y}\right)
_{1})}\varphi \leq \mu _{x}\Leftrightarrow y_{0}^{2}y_{1}^{2}\left\| \frac{1%
}{t}\right\| _{L^{1}(\left( \eta _{y}\right) _{2}^{2})}\varphi \leq \mu _{x}$
\\ 
\\ 
$\Leftrightarrow y_{0}^{2}\left\{ x^{2}\left\| \frac{1}{t}\right\|
_{L^{1}(\eta )}\widetilde{\xi }+\left( y_{1}^{2}\left\| \frac{1}{t}\right\|
_{L^{1}(\left( \eta _{y}\right) _{2}^{2})}-x^{2}r\left\| \frac{1}{t}\right\|
_{L^{1}(\eta )}\right) \delta _{0}\right\} \leq \mu _{x}$ \\ 
\\ 
$\Leftrightarrow y_{0}^{2}\left\{ x^{2}\beta _{1}^{2}\left\| \frac{1}{t}%
\right\| _{L^{1}(\eta _{1}^{2})}\widetilde{\xi }+\left( y_{1}^{2}\left\| 
\frac{1}{t}\right\| _{L^{1}(\left( \eta _{y}\right) _{2}^{2})}-x^{2}r\beta
_{1}^{2}\left\| \frac{1}{t}\right\| _{L^{1}(\eta _{1}^{2})}\right) \delta
_{0}\right\} \leq \mu _{x}$.%
\end{tabular}
\label{eq14}
\end{equation}%
By combining (\ref{eq13}) and (\ref{eq14}), we easily see that 
\begin{equation}
(T_{1},T_{2}^{2})|_{\mathcal{H}^{0}}\text{ is subnormal }\Leftrightarrow
y_{0}^{2}\left\| \frac{1}{t}\right\| _{L^{1}(\left( \eta _{y}\right)
_{1})}\varphi \leq \mu _{x}.  \label{eq15}
\end{equation}%
We thus have a characterization of the subnormality of $(T_{1},T_{2}^{2})|_{%
\mathcal{H}^{0}}$. \ From (\ref{eq11}) and (\ref{eq15}) it now follows that
the subnormality of $(T_{1},T_{2}^{2})$ implies the subnormality of $%
(T_{1},T_{2})$. \qed

It is straightforward from Definition \ref{tc} that a flat $2$-variable
weighted shift $\mathbf{T}\in \mathfrak{H}_{0}$ necessarily belongs to $%
\mathcal{TC}$. \ Thus, the following result is an easy consequence of
Theorem \ref{Thm1}.

\begin{corollary}
Let $\mathbf{T\equiv (}T_{1},T_{2})$ be a \textit{flat} $2$-variable
weighted shifts, that is, a $2$-variable weighted shift $\mathbf{T}\in 
\mathfrak{H}_{0}$ given by Figure \ref{flat}. Then we have $%
(T_{1},T_{2}^{2})\in \mathfrak{H}_{\infty }$ if and only if $%
(T_{1}^{2},T_{2})\in \mathfrak{H}_{\infty }$ if and only if $%
(T_{1},T_{2})\in \mathfrak{H}_{\infty }.$
\end{corollary}

For a flat, contractive $2$-variable weighted shift $\mathbf{T}\equiv
(T_{1},T_{2})$, we can give a concrete condition for the subnormality of $%
\mathbf{T}$. \ To do this, let $shift(\alpha _{0},\alpha _{1},\cdots )$ and $%
shift(\beta _{0},\beta _{1},\cdots )$ have Berger measures $\xi $ and $\eta $%
, respectively. \ Also, recall that for $0<\alpha <\beta $, $shift(\alpha
,\beta ,\beta ,...)$ is subnormal, with Berger measure $(1-\frac{\alpha ^{2}%
}{\beta ^{2}})\delta _{0}+\frac{\alpha ^{2}}{\beta ^{2}}\delta _{\beta ^{2}}$%
. \ To avoid trivial cases, and to ensure that each of $T_{1}$ and $T_{2}$
is a contraction, we need to assume that $ab^{n}<\Pi _{j=1}^{n}\beta _{j}$,
and we shall see in Theorem \ref{thm4} that we also need $\frac{a^{2}}{b^{2}}%
<\left\| \frac{1}{t}\right\| _{L^{1}(\eta _{1})}$, where $\eta _{1}$ is the
Berger measure of $shift(\beta _{1},\beta _{2},\beta _{3},\cdots )$. \
Finally, we know from \cite[Theorem 3.3]{CuYo2} and \cite[Section 5]{CuYo3}
that if $\mathbf{T\equiv (}T_{1},T_{2})$ is subnormal, then $\xi $ and $\eta 
$ are of the form 
\begin{equation}
\begin{tabular}{l}
$\xi =p\delta _{0}+q\delta _{1}+[1-(p+q)]\rho $ \\ 
$\eta =u\delta _{0}+v\delta _{b^{2}}+[1-(u+v)]\sigma $,%
\end{tabular}
\label{eta}
\end{equation}%
where $0<p,q,u,v<1,$ $p+q\leq 1$, $u+v\leq 1$, and $\rho ,\sigma $ are
probability measures with $\rho (\{0\}\cup \{1\})=0$, $\sigma (\{0\}\cup
\{b^{2}\})=0$. \ We then have:

\setlength{\unitlength}{1mm} \psset{unit=1mm} 
\begin{figure}[th]
\begin{center}
\begin{picture}(165,90)

\psline{->}(20,20)(85,20)
\psline(20,40)(83,40)
\psline(20,60)(83,60)
\psline(20,80)(83,80)
\psline{->}(20,20)(20,85)
\psline(40,20)(40,83)
\psline(60,20)(60,83)
\psline(80,20)(80,83)

\put(12,16){\footnotesize{$(0,0)$}}
\put(36.5,16){\footnotesize{$(1,0)$}}
\put(56.5,16){\footnotesize{$(2,0)$}}
\put(76.5,16){\footnotesize{$(3,0)$}}

\put(27,21){\footnotesize{$\alpha_{0}$}}
\put(47,21){\footnotesize{$\alpha_{1}$}}
\put(67,21){\footnotesize{$\alpha_{2}$}}
\put(81,21){\footnotesize{$\cdots$}}

\put(28,41){\footnotesize{$a$}}
\put(47,41){\footnotesize{$1$}}
\put(67,41){\footnotesize{$1$}}
\put(81,41){\footnotesize{$\cdots$}}

\put(26,63){\footnotesize{$\frac{ab}{\sqrt{\gamma_{1}(\eta_{1})}}$}}
\put(47,61){\footnotesize{$1$}}
\put(67,61){\footnotesize{$1$}}
\put(81,61){\footnotesize{$\cdots$}}

\put(27,81){\footnotesize{$\cdots$}}
\put(47,81){\footnotesize{$\cdots$}}
\put(67,81){\footnotesize{$\cdots$}}

\psline{->}(50,14)(70,14)
\put(60,10){$\rm{T}_1$}
\psline{->}(10, 50)(10,70)
\put(4,60){$\rm{T}_2$}

\put(11,39){\footnotesize{$(0,1)$}}
\put(11,59){\footnotesize{$(0,2)$}}
\put(11,79){\footnotesize{$(0,3)$}}

\put(20,28){\footnotesize{$\beta_{0}$}}
\put(20,48){\footnotesize{$\beta_{1}$}}
\put(20,68){\footnotesize{$\beta_{2}$}}
\put(21,81){\footnotesize{$\vdots$}}

\put(40,28){\footnotesize{$\frac{a\beta_{0}}{\sqrt{\gamma_{1}(\xi)}}$}}
\put(40,48){\footnotesize{$b$}}
\put(40,68){\footnotesize{$b$}}
\put(41,81){\footnotesize{$\vdots$}}

\put(60,28){\footnotesize{$\frac{a\beta_{0}}{\sqrt{\gamma_{2}(\xi)}}$}}
\put(60,48){\footnotesize{$b$}}
\put(60,68){\footnotesize{$b$}}
\put(61,81){\footnotesize{$\vdots$}}

%\\\\\\\\\

\psline{->}(130,14)(150,14)
\put(140,10){$\rm{T}_1$}
\psline{->}(97,50)(97,70)
\put(91,60){$\rm{T}_2$}

\psline{->}(100,20)(165,20)
\psline(100,40)(163,40)
\psline(100,60)(163,60)
\psline(100,80)(163,80)

\psline{->}(100,20)(100,85)
\psline(120,20)(120,83)
\psline(140,20)(140,83)
\psline(160,20)(160,83)

\put(93,16){\footnotesize{$(0,0)$}}
\put(116,16){\footnotesize{$(1,0)$}}
\put(136,16){\footnotesize{$(2,0)$}}
\put(156,16){\footnotesize{$(3,0)$}}

\put(108,21){\footnotesize{$x$}}
\put(127,21){\footnotesize{$1$}}
\put(147,21){\footnotesize{$1$}}
\put(161,21){\footnotesize{$\cdots$}}

\put(108,41){\footnotesize{$a$}}
\put(127,41){\footnotesize{$1$}}
\put(147,41){\footnotesize{$1$}}
\put(161,41){\footnotesize{$\cdots$}}

\put(108,61){\footnotesize{$a$}}
\put(127,61){\footnotesize{$1$}}
\put(147,61){\footnotesize{$1$}}
\put(161,61){\footnotesize{$\cdots$}}

\put(107,81){\footnotesize{$\cdots$}}
\put(127,81){\footnotesize{$\cdots$}}
\put(147,81){\footnotesize{$\cdots$}}
\put(161,81){\footnotesize{$\cdots$}}

\put(100,28){\footnotesize{$y$}}
\put(100,48){\footnotesize{$\beta_{1}$}}
\put(100,68){\footnotesize{$\beta_{2}$}}
\put(101,81){\footnotesize{$\vdots$}}

\put(120,28){\footnotesize{$\frac{ay}{x}$}}
\put(120,48){\footnotesize{$\beta_{1}$}}
\put(120,68){\footnotesize{$\beta_{2}$}}
\put(121,81){\footnotesize{$\vdots$}}

\put(140,28){\footnotesize{$\frac{ay}{x}$}}
\put(140,48){\footnotesize{$\beta_{1}$}}
\put(140,68){\footnotesize{$\beta_{2}$}}
\put(141,81){\footnotesize{$\vdots$}}

\end{picture}
\end{center}
\caption{Weight diagrams of the 2-variable weighted shifts in Theorem \ref%
{thm4} and Lemma \ref{exam}, respectively.}
\label{flat}
\end{figure}
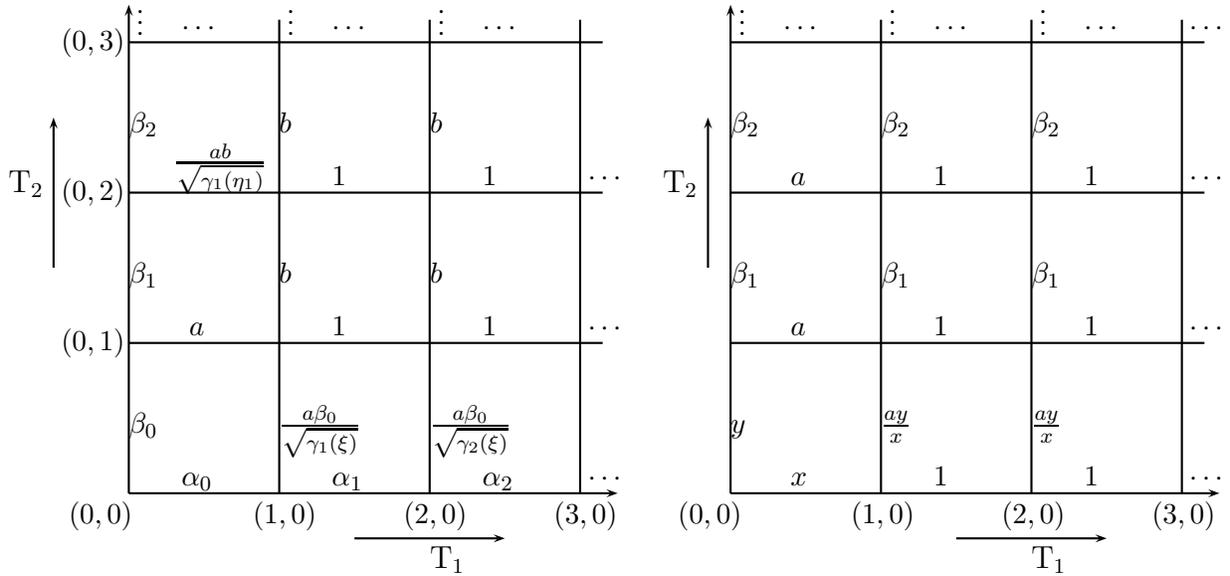

\begin{theorem}
\label{thm4}Let $\mathbf{T\equiv (}T_{1},T_{2})\in \mathfrak{H}_{0}$ be a
contractive $2$-variable weighted shift whose weight diagram is given by
Figure \ref{flat}, let $v:=\eta (\{b^{2}\})$ and $\xi \equiv p\delta
_{0}+q\delta _{1}+[1-(p+q)]\rho ,$ with $p,q>0,\;p+q\leq 1$ (cf. (\ref{eta}%
)), and let $\eta _{1}$ denote the Berger measure of $shift(\beta _{1},\beta
_{2},\cdots )$. \ Then $(T_{1},T_{2})\in \mathfrak{H}_{\infty }$ if and only
if 
\begin{equation*}
\beta _{0}\leq \min \left\{ \frac{b}{a}\sqrt{v},\sqrt{\frac{p}{(\left\| 
\frac{1}{t}\right\| _{L^{1}(\eta _{1})}-\frac{a^{2}}{b^{2}})}},\frac{b}{a}%
\sqrt{q},\sqrt{\frac{1}{\left\| \frac{1}{t}\right\| _{L^{1}(\eta _{1})}}}%
\right\} .
\end{equation*}
\end{theorem}

\textbf{Proof.} \ We first observe that%
\begin{equation}
\mu _{\mathcal{M}}=a^{2}\delta _{1}\times \delta _{b^{2}}+\delta _{0}\times
(\eta _{1}-a^{2}\delta _{b^{2}}).  \label{eq19}
\end{equation}%
Using (\ref{eta}) and (\ref{eq19}), a calculation shows that $%
(T_{1},T_{2})|_{\mathcal{M}_{1}}\in \mathfrak{H}_{\infty }$ if and only if $%
\beta _{0}\leq \frac{b}{a}\sqrt{v}$. \ Observe that 
\begin{equation*}
(\mu _{\mathcal{M}})_{ext}^{X}=\frac{1}{\left\| \frac{1}{t}\right\|
_{L^{1}(\eta _{1})}}\left\{ \left( \left\| \frac{1}{t}\right\| _{L^{1}(\eta
_{1})}-\frac{a^{2}}{b^{2}}\right) \delta _{0}+\frac{a^{2}}{b^{2}}\delta
_{1}\right\} .
\end{equation*}%
By \cite[Theorem 5.2]{CuYo1}, $(T_{1},T_{2})|_{\mathcal{N}_{1}}\in \mathfrak{%
H}_{\infty }$. \ Therefore%
\begin{equation}
\begin{tabular}{l}
$(T_{1},T_{2})\in \mathfrak{H}_{\infty }\Leftrightarrow y_{0}^{2}\left\| 
\frac{1}{t}\right\| _{L^{1}(\mu _{\mathcal{M}})}(\mu _{\mathcal{M}%
})_{ext}^{X}\leq \nu \;$(by Lemma \ref{backext}) and $\beta _{0}\leq \frac{b%
}{a}\sqrt{v}$ \\ 
\\ 
$\Leftrightarrow \beta _{0}^{2}\left\{ \left( \left\| \frac{1}{t}\right\|
_{L^{1}(\eta _{1})}-\frac{a^{2}}{b^{2}}\right) \delta _{0}+\frac{a^{2}}{b^{2}%
}\delta _{1}\right\} \leq \xi \;$and $\beta _{0}\leq \frac{b}{a}\sqrt{v}$ \\ 
\\ 
$\Leftrightarrow \beta _{0}\leq \min \left\{ \frac{b}{a}\sqrt{v},\sqrt{\frac{%
p}{(\left\| \frac{1}{t}\right\| _{L^{1}(\eta _{1})}-\frac{a^{2}}{b^{2}})}},%
\frac{b}{a}\sqrt{q},\sqrt{\frac{1}{\left\| \frac{1}{t}\right\| _{L^{1}(\eta
_{1})}}}\right\} .$%
\end{tabular}
\label{eq20}
\end{equation}%
\qed

\section{Subnormality for Powers of Hyponormal Pairs}

In this section we study the connection between the joint subnormality of
pairs $(T_{1},T_{2})\in \mathfrak{H}_{1}$ and the subnormality of the
associated monomials $T_{1}^{m}T_{2}^{n}\;(m,n\geq 1)$. \ Our results will
further exhibit the large gap between the classes $\mathfrak{H}_{\infty }$
(subnormal pairs) and $\mathfrak{H}_{0}$ (commuting pairs of subnormal
operators). \ We begin with the following proposition, which is a direct
consequence of a well known result of J. Stampfli's (\cite{Sta}, \cite{Sta2}%
): if $T$ is hyponormal and $T^{n}$ is normal for some $n\geq 1$, then $T$
is necessarily normal.

\begin{proposition}
\label{prop41}Let $\mathbf{T\equiv (}T_{1},T_{2})$ be hyponormal, and assume
that $\mathbf{(}T_{1}^{m},T_{2}^{n})$ is normal for some $m\geq 1$ and $%
n\geq 1$. Then $\mathbf{(}T_{1},T_{2})$ is normal.
\end{proposition}

In view of Proposition \ref{prop41}, one might conjecture that if $%
(T_{1},T_{2})$ is hyponormal and $T_{1}^{m}T_{2}^{n}$ is normal for some $%
m\geq 1$ and $n\geq 1$, then $(T_{1},T_{2})$ is normal (cf.\cite{Sta2}). \
But this is not true even if we assume that $(T_{1},T_{2})$ is subnormal and 
$T_{1}^{m}T_{2}^{n}$ is normal for all $m\geq 1$ and $n\geq 1$, as the
following example shows.

\begin{example}
Let $T_{1}:=U_{+}\bigoplus 0_{\infty }$ and $T_{2}:=0_{\infty }\bigoplus
U_{+}$, then $(T_{1},T_{2})$ is subnormal and $T_{1}^{m}T_{2}^{n}$ is normal
for all $m\geq 1$ and $n\geq 1$. \ However, $(T_{1},T_{2})$ is not normal.
\end{example}

Whether Proposition \ref{prop41} holds with ``normal'' replaced by
``subnormal'' is not at all obvious. \ Our main result of this section
states that it is indeed possible to have a pair $(T_{1},T_{2})\in \mathfrak{%
H}_{1}$ with $T_{1}^{m}T_{2}^{n}$ subnormal for all $m\geq 1$ and $n\geq 1$,
but such that $(T_{1},T_{2})\notin \mathfrak{H}_{\infty }$. \ (Observe,
however, that the subnormality of the monomials $T_{1}^{m}T_{2}^{n}$ is a
condition weaker than the subnormality of the pairs $\mathbf{(}%
T_{1}^{m},T_{2}^{n})$.) \ To do so, consider a subnormal weighted shift $%
shift(\beta _{1},\beta _{2},\cdots )$ with Berger measure $\eta $. \ For $%
0<a<x<1$ and $y>0$, let 
\begin{equation*}
\alpha (\mathbf{k}):=%
\begin{cases}
x & \text{if }k_{1}=0\text{ and }k_{2}=0 \\ 
a & \text{if }k_{1}=0\text{ and }k_{2}\geq 1 \\ 
1 & \text{otherwise }%
\end{cases}%
\end{equation*}%
and 
\begin{equation*}
\beta (\mathbf{k}):=%
\begin{cases}
\beta _{k_{2}} & \text{if }k_{2}\geq 1 \\ 
y & \text{if }k_{1}=0\text{ and }k_{2}=0 \\ 
\frac{ay}{x} & \text{if }k_{1}\geq 1\text{ and }k_{2}=0,%
\end{cases}%
\end{equation*}%
($\mathbf{k}=(k_{1},k_{2})\in \mathbb{Z}_{+}^{2}$). \ We now let $\mathbf{T}%
:=(T_{1},T_{2})${\normalsize \ }denote the pair of $2$-variable weighted
shift on $\ell ^{2}(\mathbb{Z}_{+}^{2})$ defined by $\alpha (\mathbf{k})$
and $\beta (\mathbf{k})$. \ We then have:

\begin{lemma}
\label{exam} Let $\mathbf{T}\equiv (T_{1},T_{2})$ be the $2$-variable
weighted shift associated with $\alpha $ and $\beta $ above. \ Then\newline
(i) $\ \mathbf{T\equiv (}T_{1},T_{2})\in \mathfrak{H}_{1}$ if and only if 
\begin{equation*}
y\leq \min \left\{ \frac{\beta _{1}x\sqrt{1-x^{2}}}{\sqrt{%
x^{2}+a^{4}-2a^{2}x^{2}}},\sqrt{\left\| \frac{1}{t}\right\| _{L^{1}(\eta
)}^{-1}}\right\} .
\end{equation*}%
\newline
(ii) \ $\mathbf{T\equiv (}T_{1},T_{2})\in \mathfrak{H}_{\infty }$ if and
only if 
\begin{equation*}
y\leq \sqrt{\left\| \frac{1}{t}\right\| _{L^{1}(\eta )}^{-1}}\cdot \sqrt{%
\frac{1-x^{2}}{1-a^{2}}}.
\end{equation*}
\end{lemma}

\textbf{Proof.} \ First observe that if $shift(y,\beta _{1},\beta
_{2},\cdots )$ is subnormal then $T_{2}$ is subnormal. \ To guarantee this,
by Lemma \ref{backext} we must have $y\leq \sqrt{\left\| \frac{1}{t}\right\|
_{L^{1}(\eta )}^{-1}}$. \ For the hyponormality of $\mathbf{(}T_{1},T_{2})$,
it suffices to apply the Six-point Test to $\mathbf{k=(}0,0)$, since 
\begin{equation*}
\mathcal{R}_{10}(\mathbf{T})\equiv \mathbf{(}T_{1},T_{2})|_{\mathcal{M}%
_{1}}\cong (I\otimes U_{+},shift(\frac{ay}{x},\beta _{1},\beta _{2},\cdots
)\otimes I)\in \mathfrak{H}_{\infty }
\end{equation*}%
and 
\begin{equation*}
\mathcal{R}_{01}(\mathbf{T})\equiv \mathbf{(}T_{1},T_{2})|_{\mathcal{N}%
_{1}}\cong (I\otimes S_{a},shift(\beta _{1},\beta _{2},\beta _{3},\cdots
)\otimes I)\in \mathfrak{H}_{\infty }.
\end{equation*}%
Thus, 
\begin{equation*}
\begin{tabular}{l}
$\left( 
\begin{array}{cc}
1-x^{2} & \frac{a^{2}y}{x}-xy \\ 
\frac{a^{2}y}{x}-xy & \beta _{1}^{2}-y^{2}%
\end{array}%
\right) \geq 0$ (by Lemma \ref{joint hypo}) \\ 
\\ 
$\Leftrightarrow y^{2}(1+\frac{a^{4}}{x^{2}}-2a^{2})\leq \beta
_{1}^{2}(1-x^{2})$ \\ 
\\ 
$\Leftrightarrow y\leq \frac{\beta _{1}x\sqrt{1-x^{2}}}{\sqrt{%
x^{2}+a^{4}-2a^{2}x^{2}}}.$%
\end{tabular}%
\end{equation*}%
Therefore, $\mathbf{T\equiv (}T_{1},T_{2})\in \mathfrak{H}_{1}$ if and only
if 
\begin{equation*}
y\leq \min \left\{ \frac{\beta _{1}x\sqrt{1-x^{2}}}{\sqrt{%
x^{2}+a^{4}-2a^{2}x^{2}}},\sqrt{\left\| \frac{1}{t}\right\| _{L^{1}(\eta
)}^{-1}}\right\} .
\end{equation*}%
We now study the subnormality of $\mathbf{T}$. $\ $Since $\mu _{\mathcal{M}%
}(s,t)=[(1-a^{2})\delta _{0}(s)+a^{2}\delta _{1}(s)]\cdot \eta (t)$ is the
Berger measure of $(I\otimes S_{a},shift(\beta _{1},\beta _{2},\beta
_{3},\cdots )\otimes I),$ Lemma \ref{backext} implies that 
\begin{equation*}
\begin{tabular}{l}
$\mathbf{T}\text{ is subnormal}$ \\ 
\\ 
$\Leftrightarrow y^{2}\left\| \frac{1}{t}\right\| _{L^{1}(\mu _{\mathcal{M}%
})}\mu _{\mathcal{M}}(s,t)_{ext}^{X}\leq (1-x^{2})\delta _{0}(s)+x^{2}\delta
_{1}(s)$ and $y\leq \sqrt{\left\| \frac{1}{t}\right\| _{L^{1}(\eta )}^{-1}}$
\\ 
\\ 
$\Leftrightarrow y^{2}\left\| \frac{1}{t}\right\| _{L^{1}(\eta
)}[(1-a^{2})\delta _{0}(s)+a^{2}\delta _{1}(s)]\leq (1-x^{2})\delta
_{0}(s)+x^{2}\delta _{1}(s)$ \\ 
\ \ \ \ \ and $y\leq \sqrt{\left\| \frac{1}{t}\right\| _{L^{1}(\eta )}^{-1}}$
\\ 
\\ 
$\Leftrightarrow y\leq \min \left\{ \sqrt{\left\| \frac{1}{t}\right\|
_{L^{1}(\eta )}^{-1}}\cdot \sqrt{\frac{1-x^{2}}{1-a^{2}}},\sqrt{\left\| 
\frac{1}{t}\right\| _{L^{1}(\eta )}^{-1}}\cdot \frac{x}{a},\sqrt{\left\| 
\frac{1}{t}\right\| _{L^{1}(\eta )}^{-1}}\right\} $ \\ 
\\ 
$\Leftrightarrow y\leq \sqrt{\left\| \frac{1}{t}\right\| _{L^{1}(\eta )}^{-1}%
}\cdot \sqrt{\frac{1-x^{2}}{1-a^{2}}}$ \\ 
\ \ \ \ \ (because $x>a$ implies $\sqrt{\frac{1-x^{2}}{1-a^{2}}}<\frac{x}{a}$
and $\sqrt{\frac{1-x^{2}}{1-a^{2}}}<1$).%
\end{tabular}%
\end{equation*}%
\qed

We now detect the hyponormality and subnormality of the powers of $\mathbf{(}%
T_{1},T_{2})$ in Lemma \ref{exam}. \ Let 
\begin{equation*}
\mathcal{H}_{(m,i)}:=\bigvee_{j=0}^{\infty }\{e_{(mj+i,k)}:m\geq 1,\ 0\leq
i\leq m-1\ {\text{and}}\ k=0,1,2,\cdots \}.
\end{equation*}%
Then $\ell ^{2}(\mathbb{Z}_{+}^{2})\equiv \bigoplus_{i=0}^{m-1}\mathcal{H}%
_{(m,i)}.$ \ Under this decomposition, we have 
\begin{equation*}
T_{1}^{m}\cong T_{1}\bigoplus (I\otimes U_{+})\bigoplus \cdots \bigoplus
(I\otimes U_{+})
\end{equation*}%
and%
\begin{equation*}
T_{2}\cong T_{2}\bigoplus (shift(\frac{ay}{x},\beta _{1},\beta _{2},\cdots
)\otimes I)\bigoplus \cdots \bigoplus (shift(\frac{ay}{x},\beta _{1},\beta
_{2},\cdots )\otimes I).
\end{equation*}%
Thus, for all $m\geq 1$ and $n\geq 1,$ 
\begin{equation*}
\mathbf{(}T_{1}^{m},T_{2}^{n})\cong \mathbf{(}T_{1},T_{2}^{n})\bigoplus
\bigoplus_{i=1}^{m-1}(C,D),
\end{equation*}%
where $C:=I\otimes U_{+}$ and $D:=(shift(\frac{ay}{x},\beta _{1},\beta
_{2},\cdots ))^{n}\otimes I.$ \ But, since $(C,D)$ is subnormal, the
hyponormality(or subnormality) of $(T_{1}^{m},T_{2}^{n})$ is equivalent to
the hyponormality (or subnormality) of $(T_{1},T_{2}^{n})$. \ Therefore, $%
(T_{1},T_{2}^{n})$ is hyponormal (or subnormal) if and only if $%
(T_{1}^{m},T_{2}^{n})$ is hyponormal (or subnormal) for all $m\geq 1.$

\begin{theorem}
\label{equivalent}For the $2$-variable weighted shift $\mathbf{T\equiv (}%
T_{1},T_{2})$ in Lemma \ref{exam}, the following are equivalent.\newline
(i) $\ T_{1}^{m}T_{2}^{n}$ is subnormal for all $m\geq 1$ and $n\geq 1$;%
\newline
(ii) \ $T_{1}T_{2}^{n}$ is subnormal for all $n\geq 1$;\newline
(iii) \ The $shift(\frac{ay\cdot \Pi _{j=1}^{n-1}\beta _{j}}{x},\Pi
_{j=n}^{2n-1}\beta _{j},\Pi _{j=2n}^{3n-1}\beta _{j},\cdots )$ is subnormal
for all $n\geq 1$;\newline
(iv) \ $y\leq \frac{x}{a}\cdot \frac{1}{\Pi _{j=1}^{n-1}\beta _{j}}\sqrt{%
\left\| \frac{1}{t}\right\| _{L^{1}(\eta ^{(n)})}^{-1}}$ for all $n\geq 1$,
where $d\eta ^{(n)}(t):=\frac{t^{1-\frac{1}{n}}}{\beta _{1}^{2}\cdots \beta
_{n-1}^{2}}d\eta (t^{\frac{1}{n}})$.
\end{theorem}

\textbf{Proof.} $\ (i)\Longleftrightarrow (ii)$ \ From the above
observations, we can see that $T_{1}^{m}T_{2}^{n}$ is subnormal for all $%
m\geq 1$ and $n\geq 1$ if and only if $T_{1}T_{2}^{n}$ and $CD$ are
subnormal for all $n\geq 1$. $\ $But observe that $CD$ is always subnormal
if $shift(\frac{ay}{x},\beta _{1},\beta _{2},\cdots )$ is subnormal.

$(ii)\Longleftrightarrow (iii)$ \ Let $\mathcal{M}_{(i,j)}:=\bigvee
\{e_{i+k,j+k}:k=0,1,2,\cdots \}$ for $i,j\geq 0$ with $ij=0.$ \ Then $\ell
^{2}(\mathbb{Z}_{+}^{2})\equiv \bigoplus_{i,j=0}^{\infty }\mathcal{M}%
_{(i,j)} $. $\ $Under this decomposition, we have 
\begin{equation*}
T_{1}T_{2}^{n}\cong \cdots \bigoplus W_{-1}\bigoplus W_{0}\bigoplus
W_{1}\bigoplus \cdots ,
\end{equation*}%
where 
\begin{eqnarray*}
W_{-1} &:&=shift(a\Pi _{j=n}^{2n-1}\beta _{j},\Pi _{j=2n}^{3n-1}\beta
_{j},\Pi _{j=3n}^{4n-1}\beta _{j},\cdots ):\mathcal{M}_{(0,1)}%
\longrightarrow \mathcal{M}_{(0,1)}, \\
W_{0} &:&=shift(ay\cdot \Pi _{j=1}^{n-1}\beta _{j},\Pi _{j=n}^{2n-1}\beta
_{j},\Pi _{j=2n}^{3n-1}\beta _{j},\cdots ):\mathcal{M}_{(0,0)}%
\longrightarrow \mathcal{M}_{(0,0)}\text{, and} \\
W_{1} &:&=shift(\frac{ay}{x}\cdot \Pi _{j=1}^{n-1}\beta _{j},\Pi
_{j=n}^{2n-1}\beta _{j},\Pi _{j=2n}^{3n-1}\beta _{j},\cdots ):\mathcal{M}%
_{(1,0)}\longrightarrow \mathcal{M}_{(1,0)}.
\end{eqnarray*}%
Since $W_{-1}$ is subnormal, the result follows from the fact that if $W_{1}$
is subnormal then $W_{0}$ is also subnormal.

$(iii)\Longleftrightarrow (iv)$ \ Since $shift(\beta _{1},\beta _{2},\beta
_{3},\cdots )$ has Berger measure $\eta $, we can use mathematical induction
to show that $shift(\beta _{n},\beta _{n+1},\beta _{n+2},\cdots )$ has
Berger measure $\frac{t^{n-1}}{\beta _{1}^{2}\cdots \beta _{n-1}^{2}}d\eta
(t)$ for each $n\geq 1$. Thus by Lemma \ref{multiply}, $shift(\Pi
_{j=n}^{2n-1}\beta _{j},\Pi _{j=2n}^{3n-1}\beta _{j},\Pi _{j=3n}^{4n-1}\beta
_{j},\cdots )$ has Berger measure $d\eta ^{(n)}(t)\equiv \frac{t^{1-\frac{1}{%
n}}}{\beta _{1}^{2}\cdots \beta _{n-1}^{2}}d\eta (t^{\frac{1}{n}})$ for each 
$n\geq 1.$ \ Therefore, by Lemma \ref{backward} we see that \newline
$shift(\frac{ay\cdot \Pi _{j=1}^{n-1}\beta _{j}}{x},\Pi _{j=n}^{2n-1}\beta
_{j},\Pi _{j=2n}^{3n-1}\beta _{j},\cdots )$ is subnormal if and only if $%
y\leq \frac{x}{a}\cdot \frac{1}{\Pi _{j=1}^{n-1}\beta _{j}}\sqrt{\left\| 
\frac{1}{t}\right\| _{L^{1}(\eta ^{(n)})}^{-1}}$. \qed

For a concrete example, let $d\eta (t):=dt$ on $[\frac{1}{2},\frac{3}{2}]$,
so that $\beta _{1}=1$ and $\left\| \frac{1}{t}\right\| _{L^{1}(\eta )}=\ln
3 $. \ Since $\gamma _{n-1}=\beta _{1}^{2}\beta _{2}^{2}\cdots \beta
_{n-1}^{2}=\int_{\frac{1}{2}}^{\frac{3}{2}}t^{n-1}d\eta (t)=\frac{1}{n}(%
\frac{3^{n}-1}{2^{n}})$ and $\gamma _{2n-1}=\frac{1}{2n}(\frac{3^{2n}-1}{%
2^{2n}})$, it follows that $shift(\beta _{n},\beta _{n+1},\cdots )$ has
Berger measure $\frac{n\cdot 2^{n}\cdot t^{n-1}}{3^{n}-1}dt$ for each $n\geq
1$ on $[\frac{1}{2},\frac{3}{2}]$ and \newline
$shift(\Pi _{j=n}^{2n-1}\beta _{j},\Pi _{j=2n}^{3n-1}\beta _{j},\cdots )$
has Berger measure $d\eta ^{(n)}(t)=\frac{2^{n}}{3^{n}-1}dt$ on $[\left( 
\frac{1}{2}\right) ^{n},\left( \frac{3}{2}\right) ^{n}]\;($all $n\geq 1)$. \
Moreover, $\sqrt{\left\| \frac{1}{t}\right\| _{L^{1}(\eta ^{(n)})}^{-1}}=%
\sqrt{\frac{3^{n}-1}{n2^{n}\ln 3}}$. \ Thus, Lemma \ref{exam} implies that%
\newline
(i) \ $T_{1}$ is subnormal if $0<a<x<1$;\newline
(ii) \ $T_{2}$ is subnormal $\Leftrightarrow $ $y\leq \sqrt{\frac{1}{\ln 3}}$%
;\newline
(iii) \ $\mathbf{(}T_{1},T_{2})\in \mathfrak{H}_{1}\Leftrightarrow y\leq
m:=\min \{\frac{x\sqrt{1-x^{2}}}{\sqrt{x^{2}+a^{4}-2a^{2}x^{2}}},\sqrt{\frac{%
1}{\ln 3}}\}$;\newline
(iv) $\ \mathbf{(}T_{1},T_{2})\in \mathfrak{H}_{\infty }\Leftrightarrow
y\leq s:=\sqrt{\frac{1}{\ln 3}\frac{1-x^{2}}{1-a^{2}}}$.$\newline
$Therefore, we have the following result.

\begin{example}
\label{four}For $s<y\leq m$ and $0<a<x<1$, we have\newline
(i) \ $\mathbf{T}\equiv (T_{1},T_{2})\in \mathfrak{H}_{1}$;\newline
(ii) \ $\mathbf{T}\equiv (T_{1},T_{2})\notin \mathfrak{H}_{\infty }$;\newline
(iii) \ $T_{1}^{m}T_{2}^{n}$ is subnormal for all $m\geq 1,n\geq 1.$\newline
For, observe that if $0<a<x<1$, then $s\equiv \sqrt{\frac{1}{\ln 3}\frac{%
1-x^{2}}{1-a^{2}}}<\frac{x\sqrt{1-x^{2}}}{\sqrt{x^{2}+a^{4}-2a^{2}x^{2}}}$
and $s<\sqrt{\frac{1}{\ln 3}}$; thus, $s<m$, and it is then possible to
choose values of $y$ between these two quantities. \ From Theorem \ref%
{equivalent}, we can see that $T_{1}^{m}T_{2}^{n}$ is subnormal for all $%
m\geq 1,n\geq 1$ if and only if $y\leq \frac{x}{a}\cdot \frac{1}{\Pi
_{j=1}^{n-1}\beta _{j}}\sqrt{\left\| \frac{1}{t}\right\| _{L^{1}(\mu _{\eta
})}^{-1}}=\frac{x}{a}\sqrt{\frac{1}{\ln 3}}.$ \ But since $y\leq \sqrt{\frac{%
1}{\ln 3}}<\frac{x}{a}\sqrt{\frac{1}{\ln 3}}$, it follows that $%
T_{1}^{m}T_{2}^{n}$ is subnormal for all $m\geq 1,n\geq 1.$
\end{example}


\begin{thebibliography}{99}
\bibitem{Abr} M. Abrahamse, Commuting subnormal operators, \textit{Ill.
Math. J.} 22(1978), 171-176.

\bibitem{Ath} A. Athavale, On joint hyponormality of operators, \textit{\
Proc. Amer. Math. Soc}. 103(1988), 417-423.

\bibitem{Atk} K. Atkinson, \textit{Introduction to Numerical Analysis},
Wiley and Sons, 2nd. Ed. \ 1989.

\bibitem{Bra} J. Bram, Subnormal operators, \textit{Duke Math. J.} 2(1955),
75-94.

\bibitem{Con} J. Conway, \textit{The Theory of Subnormal Operators,}
Mathematical Surveys and Monographs, vol. 36, Amer. Math. Soc., Providence,
1991.

\bibitem{bridge} R. Curto, Joint hyponormality: A bridge between
hyponormality and subnormality, \textit{Proc. Symposia Pure Math.} 51(1990),
69-91.

\bibitem{QHWS} R. Curto, Quadratically hyponormal weighted shifts, \textit{%
Integral Equations Operator Theory} 13(1990), 49-66.

\bibitem{OTAMP} R. Curto, An operator-theoretic approach to truncated moment
problems, in \textit{Linear Operators}, Banach Center Publ., vol. 38, 1997,
pp. 75-104.

\bibitem{CLL2} R. Curto, S.H. Lee and W.Y. Lee, A new criterion for $k$%
-hyponormality via weak subnormality, \textit{Proc. Amer. Math. Soc.}
133(2005), 1805-1816.

\bibitem{CLY} R. Curto, S.H. Lee and J. Yoon, $k$-hyponormality of
multivariable weighted shifts, \textit{J. Funct. Anal.} 229(2005), 462-480.

\bibitem{CuLe1} R. Curto and W.Y. Lee, Joint hyponormality of Toeplitz
pairs, \textit{Memoirs Amer. Math. Soc.} no. 712, Amer. Math. Soc.,
Providence, 2001.

\bibitem{CuLe2} R. Curto and W.Y. Lee, Towards a model theory for $2$%
-hyponormal operators, \textit{Integral Equations Operator Theory} 44(2002),
290-315.

\bibitem{CMX} R. Curto, P. Muhly and J. Xia, Hyponormal pairs of commuting
operators, \textit{Operator Theory: Adv. Appl.} 35(1988), 1-22.

\bibitem{CuP} R. Curto and S. Park, $k$-hyponormality of powers of weighted
shifts, \textit{Proc. Amer. Math. Soc.} 131(2002), 2762-2769.

\bibitem{CPY} R. Curto, Y. Poon and J. Yoon, Subnormality of Bergman-like
weighted shifts, \textit{J. Math. Anal. Appl.} 308(2005), 334-342.

\bibitem{CuYo1} R. Curto and J. Yoon, Jointly hyponormal pairs of subnormal
operators need not be jointly subnormal, \textit{Trans. Amer. Math. Soc.}
358(2006), 5139-5159.

\bibitem{CuYo2} R. Curto and J. Yoon, Disintegration-of-measure techniques
for multivariable weighted shifts, \textit{Proc. London Math. Soc.}
93(2006), 381-402.

\bibitem{CuYo3} R. Curto and J. Yoon, Propagation phenomena for hyponormal $%
2 $-variable weighted shifts, \textit{J. Operator Theory}, to appear.

\bibitem{Fra} E. Franks, Polynomially subnormal operator tuples, \textit{J.
Operator Theory} 31(1994), 219-228.

\bibitem{GeWa} Gellar and Wallen, Subnormal weighted shifts and the
Halmos-Bram criterion,\textit{\ Proc. Japan Acad.} 46(1970), 375-378.

\bibitem{Hal} P.R. Halmos, \textit{A Hilbert Space Problem Book},
Springer-Verlag, Berlin and New York.

\bibitem{JeLu} N.P. Jewell and A.R. Lubin, Commuting weighted shifts and
analytic function theory in several variables, \textit{J. Operator Theory}
1(1979), 207-223.

\bibitem{Lu1} A. Lubin, Weighted shifts and product of subnormal operators, 
\textit{Indiana Univ. Math. J.} 26(1977), 839-845.

\bibitem{Lu2} A. Lubin, Extensions of commuting subnormal operators, \textit{%
Lecture Notes in Math.} 693(1978), 115-120.

\bibitem{Lu3} A. Lubin, A subnormal semigroup without normal extension, 
\textit{Proc. Amer. Math. Soc.} 68(1978), 176-178.

\bibitem{Pau} V. Paulsen, Completely bounded maps and dilations, \textit{%
Pitmam Research Notes in Mathematics Series}, vol. 146, Longman Sci. Tech.,
New York, 1986.

\bibitem{Smu} J.L. Smul'jan, An operator Hellinger integral, \textit{Mat.
Sb. (N.S.)} 49 (1959), 381--430 (in Russian).

\bibitem{Sta} J. Stampfli, Hyponormal operators, \textit{Pacific J. Math.}
12(1962), 1453-1458.

\bibitem{Sta2} J. Stampfli, Which weighted shifts are subnormal?, \textit{%
Pacific J. Math}. 17(1966), 367-379.

\bibitem{Wol} Wolfram Research, Inc. \textit{Mathematica}, Version 4.2, 
\textit{Wolfram Research Inc.}, Champaign, IL, 2002.

\bibitem{Yoo} J. Yoon, Disintegration of measures and contractive $2$%
-variable weighted shifts, preprint 2006.

\end{thebibliography}
\end{document}